\documentclass[12pt]{article}

\usepackage{amsmath,amsthm,amsfonts,amssymb}

\newtheorem{theorem}{Theorem}[section]

\newtheorem{proposition}[theorem]{Proposition}
\newtheorem{lemma}[theorem]{Lemma}
\newtheorem{definition}[theorem]{Definition}
\newtheorem{remark}[theorem]{Remark}

\def\cB{\mathcal{B}}

\def\cD{\mathcal{D}}
\def\cE{\mathcal{E}}

\def\cH{\mathcal{H}}
\def\cF{\mathcal{F}}

\def\cM{\mathcal{M}}
\def\cN{\mathcal{N}}

\def\cP{\mathcal{P}}
\def\cS{\mathcal{S}}
\def\cX{\mathcal{X}}

\def\bC{\mathbb{C}}

\def\bR{\mathbb{R}}

\def\intR{\int_{\mathbb R}}
\def\inte#1{\int_0^#1}

\newcommand*\samethanks[1][\value{footnote}]{\footnotemark[#1]}
\begin{document}

\title{SPDEs with fractional noise in space\\
with index $H<1/2$}

\author{Raluca Balan\footnote{Department of Mathematics and Statistics, University of Ottawa,
585 King Edward Avenue, Ottawa, ON, K1N 6N5, Canada. E-mail
address: rbalan@uottawa.ca. Research supported by a
grant from the Natural Sciences and Engineering Research Council
of Canada.}\and
Maria Jolis\thanks{Departament de Matem\`atiques,
Universitat Aut\`{o}noma de Barcelona, 08193 Bellaterra (Barcelona), Catalonia, Spain.
E-mail addresses: mjolis@mat.uab.cat, quer@mat.uab.cat. Research supported by the grant
MCI-FEDER MTM2012-33937. Corresponding author: Llu\'{i}s Quer-Sardanyons.}\and
Llu\'{i}s Quer-Sardanyons \samethanks
}

\date{July 15, 2014}
\maketitle

\begin{abstract}
\noindent In this article, we consider the stochastic wave and heat equations on $\bR$
with non-vanishing initial conditions, driven by a Gaussian noise which is white in
time and behaves in space like a fractional Brownian motion of index $H$, with $1/4<H<1/2$.
We assume that the diffusion coefficient is given by an affine function $\sigma(x)=ax+b$,
and the initial value functions are bounded and H\"older continuous of order $H$.
We prove the existence and uniqueness of the mild solution for both equations. We show that the solution is $L^{2}(\Omega)$-continuous and its $p$-th moments are uniformly bounded, for any $p \geq 2$.
\end{abstract}

\noindent {\em MSC 2010:} Primary 60H15; secondary 60H05

\vspace{3mm}

\noindent {\em Keywords:} stochastic wave equation; stochastic heat equation; fractional Brownian motion; random field solution 
\section{Introduction}

In this article, we consider the stochastic wave equation:
\begin{equation}
\left\{\begin{array}{rcl}
\displaystyle \frac{\partial^2 u}{\partial t^2}(t,x) & = & \displaystyle \frac{\partial^2 u}{\partial x^2}(t,x)+\sigma(u(t,x))\dot{X}(t,x), \quad t\in [0,T], \ x \in \bR \\[2ex]
\displaystyle u(0,x) & = & u_0(x), \\[1ex]
\displaystyle \frac{\partial u}{\partial t}(0,x) & = & v_0(x),
\end{array}\right. \label{wave} \tag{SWE}
\end{equation}
and the stochastic heat equation:
\begin{equation}
\left\{\begin{array}{rcl}
\displaystyle \frac{\partial u}{\partial t}(t,x) & = &\displaystyle  \frac{1}{2}\, \frac{\partial^2 u}{\partial x^2}(t,x) + \sigma(u(t,x))\dot{X}(t,x), \quad t\in [0,T], \ x \in \bR \\[2ex]
\displaystyle u(0,x) & = & u_0(x)
\end{array}\right. \label{heat} \tag{SHE}
\end{equation}

\noindent where $\sigma(x)=ax+b$ is an affine function and $\dot{X}$ denotes the formal derivative of a spatially homogeneous
Gaussian noise $X$, which is white in time and behaves in space like a fractional Brownian motion (fBm)
with index $H \in (1/4,1/2)$. The precise definition of $X$ is given in Section \ref{section-noise} below.
The initial value functions $u_0$ and $v_0$ are bounded and H\"older continuous of order $H$.

We denote by $G_t(x)$ the fundamental solution of the wave (respectively heat) equation, that is
$$G_t(x)=\frac{1}{2}1_{\{|x|<t\}} \quad \mbox{for the wave equation},$$
$$G_t(x)=\frac{1}{(2\pi t)^{1/2}} \exp\left(-\frac{|x|^2}{2t}\right) \quad \mbox{for the heat equation.}$$

The goal of the present article is to prove the following result.

\begin{theorem}
\label{main-thm} Equation \eqref{wave} (respectively \eqref{heat})
has a unique solution $u=\{u(t,x);\linebreak t \in [0,T], x\in \bR \}$, which
is $L^2(\Omega)$-continuous and satisfies, for any $p
\geq 2$,
$$\sup_{(t,x) \in [0,T] \times \bR}E|u(t,x)|^p<\infty$$
and
\begin{equation}
 \sup_{(t,x) \in [0,T] \times \bR} \int_0^t \int_{\bR^2}G_{t-s}^{2}(x-y)\dfrac{\Big(E|u(s,y)-u(s,z)|^p\Big)^{2/p}}{|y-z|^{2-2H}}\,dy\,dz\,ds <\infty.
 \label{eq:48}
\end{equation}
\end{theorem}

In particular, Theorem \ref{main-thm} covers the case of equation \eqref{wave}
with $\sigma(x)=x$, $u_0(x)=b$ and $v_0=0$, and equation \eqref{heat} with $\sigma(x)=x$ and $u_0(x)=b$,
which are known in the literature as the Hyperbolic
Anderson Model (HAM), respectively the Parabolic Anderson Model (PAM).
In fact, the solution $\bar{u}$ of (PAM) can be written as $\bar{u}=u+b$, where $u$ is the solution to \eqref{heat}
with $\sigma(x)=x+b$ and $u_0=0$. Equation (PAM) plays a major role in the study of the KPZ equation in physics, via the Hopf-Cole transformation. Its discrete form was studied in \cite{carmona-mol}. One possible method for studying equations (HAM) and (PAM) is based on the idea that the solution can be expressed as a series of multiple stochastic integrals with respect to $X$. This method was used in references
\cite{hu01,hu-nualart,HNS11, BT10} in the case of the heat equation, and in references \cite{DMT08,dalang-mueller09,B12} in the case of the wave equation. This approach is particularly useful when the noise behaves in time like a fBm, and martingale techniques cannot be applied.
We do not pursue this approach here. Instead, we will use the classical method of Picard iterations, our main efforts being dedicated to showing that the Picard iteration sequence is well-defined and converges (in a certain space).

The concept of solution is defined as follows. We say that a random field $u=\{u(t,x); t \in [0,T],
x\in \bR \}$ is a (mild) {\em solution} of \eqref{wave} (respectively \eqref{heat}),
if $u$ is predictable and for any $(t,x) \in [0,T] \in \bR$
\begin{equation}
\label{int-eq}
u(t,x)=w(t,x)+\int_0^t\int_{\bR} G_{t-s}(x-y)\,\sigma(u(s,y))\,X(ds,dy) \quad {\rm a.s.}
\end{equation}
where the stochastic integral is interpreted in the sense explained in Section \ref{section-stoch-int} below,
and $w=\{w(t,x); t \in [0,T],x \in \bR\}$ is the solution of the homogeneous wave (respectively heat) equation with the same initial conditions as in \eqref{wave} (respectively \eqref{heat}), namely:
$$w(t,x)=\frac{1}{2}\int_{x-t}^{x+t}v_0(y)dy+\frac{1}{2}\Big(u_0(x+t)+u_0(x-t)\Big) \quad \mbox{for the wave equation},$$
$$w(t,x)= \int_{\bR}G_t(x-y)u_0(y)dy \quad \mbox{for the heat equation}.$$

This problem has a very rich history, since stochastic partial differential equations
(SPDEs) driven by a spatially homogeneous Gaussian noise have been studied intensively
in the past fifteen years.
We recall that a spatially homogeneous Gaussian noise is a zero-mean Gaussian process $X=\{X_t(\varphi);t \geq 0,\varphi \in \cD(\bR^d)\}$ with covariance
$$E[X_t(\varphi)X_s(\psi)]=(t \wedge s)\,\Gamma(\varphi * \tilde{\psi}),$$
where $\cD(\bR^d)$ is the set of infinitely differentiable functions on $\bR^d$ with compact support, $\Gamma$ is a non-negative-definite tempered distribution on $\bR^d$ and $\tilde{\psi}(x)=\psi(-x)$. By the Bochner-Schwartz theorem, there exists a tempered measure $\mu$ on $\bR^d$ whose Fourier transform in $\cS'(\bR)$ is $\Gamma$. Here we denote by $\cS'(\bR)$ the space of tempered distributions on $\bR$.
Therefore,
\begin{equation}
\label{cov-X-noise}
E[X_t(\varphi)X_s(\psi)]=(t \wedge s)\,\int_{\bR^d}\cF \varphi(\xi) \overline{\cF \psi(\xi)}\mu(d\xi).
\end{equation}
Since its covariance is invariant under translations, the noise can be viewed as a stationary random distribution, in the sense introduced by It\^o in \cite{ito54} (see Section \ref{section-noise} below for the precise definition).

Due to difficulties in the construction of the stochastic integral with respect to $X$, most of the results
related to the study of SPDEs with this type of noise were obtained under the following assumption:

\vspace{2mm}
\noindent {\bf Assumption A.} {\em $\Gamma$ is given by a non-negative-definite tempered measure (or in particular, $\Gamma$ is given by a non-negative locally integrable function $f$).}

\vspace{2mm}

In the presence of this assumption, is is known that a general class of SPDEs with non-vanishing initial
conditions (which includes the wave equation in dimension $d \leq 3$ and the heat equation in any
dimension) have random field solutions
(see, e.g., \cite{walsh86,carmonanualart,dalang-frangos98,millet-sanz99,dalang99,dalang-quer11}.
Various properties of the solution, like H\"older continuity of the sample paths or smoothness of
the law, have been investigated by many authors
(see, e.g., \cite{bally,mms,sanz-sarra00,quer-sanz,nualart-quer07,sanz,dalang-sanz}).
Concerning the spectral measure $\mu$, in all above-mentioned references it is assumed that
\begin{equation}
 \int_{\bR^d}\frac{1}{1+|\xi|^2}\mu(d\xi)<\infty.
 \label{eq:50}
\end{equation}
On the other hand, as far as the semigroup approach to SPDEs is concerned,
we remark that Peszat and Zabczyk \cite{PZ07} (see also \cite{pz97})
obtained the existence and uniqueness of a function-space valued solution to the stochastic wave equation (with $d\leq 3$)
and stochastic heat equation (in any dimension) under condition \eqref{eq:50} and
the following:

\smallskip
\noindent {\bf Assumption B.} {\em There exists a constant $C>0$ such that $\Gamma+C\lambda_d$ is a non-negative measure, where $\lambda_d$ is the Lebesgue measure on $\bR^d$.}
\smallskip

In the case of the stochastic wave equation in any space dimension, the existence of the solution has been studied in \cite{peszat02,dalang-mueller03,conus-dalang} using different approaches. More precisely, in \cite{peszat02} the covariance $\Gamma$ is assumed to satisfy Assumption B, while in \cite{dalang-mueller03,conus-dalang}
the authors suppose that it fulfils Assumption A. On the other hand, in \cite{dalang-mueller03} the spectral measure $\mu$ satisfies \eqref{eq:50}, while in
\cite{peszat02,conus-dalang} it satisfies
\begin{equation}
 \sup_{\eta \in \bR}\int_{\bR^d}\frac{1}{1+|\xi-\eta|^2}\mu(d\xi)<\infty.
 \label{eq:51}
\end{equation}
Using arguments from \cite{karczewska-zabczyk00}, Peszat \cite{peszat02} showed that conditions \eqref{eq:50} and \eqref{eq:51} are equivalent, if Assumption B holds.

In this article, we will assume that $d=1$ and $X$ is a spatially homogeneous Gaussian noise with the same spectral measure $\mu$ as the fBm of index $H$, i.e.
\begin{equation}
\label{def-mu}
\mu(d\xi)=c_H |\xi|^{1-2H}d\xi,
\end{equation}
with
\begin{equation}
\label{def-cH}
c_H=\frac{\Gamma(2H+1) \sin(\pi H)}{2\pi}.
\end{equation}

We recall that the fBm with index $H \in (0,1)$ is a zero-mean Gaussian process $B=\{B(x)\}_{x  \in \bR}$ with covariance: (see e.g. Section 7.2.2 of \cite{ST94})
\begin{equation}
\label{cov-fBm1}
E[B(x)B(y)]=\int_{\bR}\cF 1_{[0,x]}(\xi) \overline{\cF 1_{[0,y]}(\xi)}\mu(d\xi)
\end{equation}
where $\cF$ denotes the Fourier transform, and $\mu$ is given by \eqref{def-mu}. The fBm with index $H=1/2$ coincides with the Brownian motion.

Note that the measure $\mu$ given by \eqref{def-mu} satisfies \eqref{eq:50}, for any $H \in (0,1)$.
However, condition \eqref{eq:51} {\em does not hold} when $H<1/2$ (see Appendix \ref{appendix-peszat-cond}).
On the other hand, when $H>1/2$, the Fourier transform in $\cS'(\bR)$ of $\mu$ is the locally integrable function $f(x)=H(2H-1)|x|^{2H-2}$, which satisfies Assumption A above. But when $H<1/2$, Assumption A fails, as the Fourier transform of $\mu$ in $\cS'(\bR)$ is a genuine distribution $\Gamma$, which is obtained by regularization: (see e.g. Chapter 1, Section 3 of \cite{gelfand-shilov64})
$$\Gamma(\varphi)=H(2H-1)\int_{\bR}(\varphi(x)-\varphi(0))|x|^{2H-2}dx, \quad \varphi\in \cD(\bR),$$
and coincides with $(1/2)V''$, where $V''$ denotes the second distributional derivative of $V(x)=|x|^{2H}$ (see \cite{jolis10}). Therefore, the techniques used in the references mentioned above cannot be applied in the case $H<1/2$.

The first step in the study of SPDEs is to develop a stochastic integral with respect to the noise. Since the trajectories of the fBm are $\alpha$-H\"older continuous with $\alpha<H$, the fBm with index $H<1/2$ has ``rougher'' sample paths than the Brownian motion. For this reason, we expect more restrictive conditions for integration with respect to a Gaussian noise which behaves in space like a fBm with $H<1/2$, compared to the ``smoother'' case $H>1/2$.

It was shown in \cite{jolis10} that the domain of the Wiener integral with respect to the fBm of index $H \in (0,1)$
is the completion of $\cD(\bR)$ with respect to the inner product
\begin{equation}
\label{cov-fBm2}
\langle \varphi, \psi \rangle_{\Lambda}=\int_{\bR}\cF \varphi(\xi)\overline{\cF \psi(\xi)}\mu(d\xi),
\end{equation}
and coincides with the space of distributions $S\in \cS'(\bR)$, whose Fourier transform $\cF S$ is a locally
integrable function which satisfies $\int_{\bR}|\cF S(\xi)|^2\mu(d\xi)<\infty$.
The authors of \cite{BGP12} have recently proved that a similar characterization can be given for the class of
stochastic integrands with respect to the space-time
Gaussian noise $X$. Based on this characterization, we give a new criterion for integrability with respect to $X$,
which is inherited from the theory of fractional Sobolev
spaces, and constitutes the starting point in the developments in the present article. In particular, this criterion
supplies us with the necessary tools for proving that the Picard
iteration sequence is well-defined and converges to the solution of equation \eqref{wave} (or \eqref{heat}). Furthermore, we remark that, because of the above-mentioned techniques,
the term in \eqref{eq:48} comes into the picture in a quite natural way when setting up the Picard scheme, and indeed has been crucial in order to prove the uniqueness of the solution.

The restriction $H>1/4$ arises from a technical condition that we need to impose on the fundamental solution $G$, namely: (see Remark \ref{remark-H} below)
$$\int_0^T \int_{\bR}|\cF G_t(\xi)|^2 \,|\xi|^{2(1-2H)}d\xi\,dt<\infty.$$
An intuitive explanation is that if $H \leq 1/4$, the noise is so ``rough'' in space that the Sobolev space methods of the present article cannot be applied.

This article is organized as follows. In Section \ref{section-noise-int}, we collect all the preliminary
results about the noise and the stochastic integral, together with the new criterion for integrability
mentioned above. The proof of Theorem \ref{main-thm} is presented in Section \ref{section-proof}.
Each of these sections is divided into several sub-sections, which are summarized at the beginning of the section. Some auxiliary results are presented in Appendices A, B and C.

We conclude the introduction with few words about the notation. We denote by $\cD(K)$ the space of infinitely differentiable functions on $\bR^d$ with compact support contained
in an open set $K$. We denote by $\cS(\bR)$ the class of rapidly decreasing infinitely differentiable functions on $\bR$.
The Fourier transform of a function $\varphi \in L^1(\bR^d)$ is defined by $$\cF \varphi(\xi)=\int_{\bR^d}e^{-i \xi \cdot x}\varphi(x)dx,$$
where $\xi \cdot x=\sum_{i=1}^{d}\xi_i x_i$ is the Euclidean inner product in $\bR^d$. We let $\cB(\bR^d)$ be the class of Borel sets in $\bR^d$ and $\cB_{b}(\bR^d)$ be the class of bounded Borel sets in $\bR^d$.
We denote by $L_{\bC}^2(\bR^d,\mu)$ the space of complex-valued functions on $\bR^d$, which are square-integrable with respect to the measure $\mu$.
We denote by $L^2_{\bC}(\Omega)$ the space of complex-valued square-integrable random variables defined on $\Omega$. The same notations without the subscript $\bC$ denote the corresponding subspaces consisting of real-valued elements.

\section{The noise and the stochastic integral}
\label{section-noise-int}

This section is divided in three parts.
In Section \ref{section-noise},
we introduce the Gaussian noise $X$, which can be viewed as a stationary random distribution,
as in \cite{B12-Fourier}. In Section \ref{section-stoch-int}, we recall the construction of
the stochastic integral with respect to $X$, following closely the approach of \cite{BGP12}.
We point out that the considerations in Sections \ref{section-noise} and
\ref{section-stoch-int} are indeed valid for any $H\in (0,1)$
and in fact for any symmetric measure $\mu$ on $\bR$ which satisfies \eqref{eq:50}.
In this case, the process $\{X(t,x)\}_{x \in \bR}$ (defined in Remark \ref{random-field} below) is
a Gaussisan process with stationary increments and spectral measure $t\mu$.
Finally, in Section \ref{section-new-char} we consider the case
$H\in (0,\frac12)$ and we obtain a new criterion for integrability with
respect to $X$, using tools from the theory of fractional Sobolev spaces, borrowed from \cite{DPV}.

\subsection{The noise}
\label{section-noise}

We let $X=\{X(\varphi); \varphi \in \cD((0,\infty) \times \bR)\}$ be a zero-mean Gaussian process, defined on a probability space $(\Omega,\cF,P)$, with covariance:
\begin{equation}
\label{cov-noise}
E[X(\varphi)X(\psi)]=\int_{0}^{\infty}\int_{\bR} \cF \varphi(t,\cdot)(\xi) \overline{\cF \psi(t,\cdot)(\xi)}\mu(d\xi) dt =:\langle \varphi,\psi \rangle_{\cH},
\end{equation}
for any $\varphi,\psi \in \cD((0,\infty) \times \bR)$, where $\mu$ is given by
$\mu(d\xi)=c_H\, |\xi|^{1-2H} d\xi$ (see \eqref{def-mu}) with $0<H<1$.
Note that in the space variable, $X$ has the same covariance structure as the fBm with index $H$ (see \eqref{cov-fBm2}).

\vspace{2mm}

We denote by $\cH$ the completion of $\cD((0,\infty)  \times \bR)$ with respect to $\langle \cdot, \cdot \rangle_{\cH}$. The map $\varphi \mapsto X(\varphi)$ is an isometry from $\cD((0,\infty)  \times \bR)$ to $L^2(\Omega)$ which can be extended to $\cH$. We use the notation:
$$X(\varphi)=\int_{0}^{\infty} \int_{\bR}\varphi(t,x)X(dt,dx), \quad \varphi \in \cH.$$

The process $X=\{X(\varphi); \varphi \in \cD((0,\infty)  \times \bR)\}$ is a (real) {\em stationary random distribution}, a concept which was introduced by It\^o in dimension one (see \cite{ito54}), and generalized to higher dimensions by Yaglom (see \cite{yaglom57}). This means that the map $\varphi \to X(\varphi)$ is linear and continuous from $\cD((0,\infty)  \times \bR)$ to $L^2(\Omega)$, and the covariance of $X$ is invariant under translations, i.e.
$$E[X(\tau_{h}\varphi) X(\tau_h\varphi)]=E[X(\varphi)X(\psi)] \quad \mbox{for any} \ h \in \bR_{+} \times \bR,$$
where $(\tau_h \varphi)(t,x)=\varphi(t+h_1,x+h_2)$ and $h=(h_1,h_2)$.
By Theorem 3 of \cite{yaglom57}, $X$ has the {\em spectral representation}:
$$X(\varphi)=\int_{\bR}\int_{\bR}\cF\varphi(\tau,\xi) \cM(d\tau,d\xi), \quad \varphi \in \cD((0,\infty) \times \bR),$$
where $\cM=\{\cM(A);A \in \cB_{b}(\bR^2)\}$ is a Gaussian complex random measure on $\bR^2$ with zero-mean and control measure $\Pi(d\tau,d\xi)=\frac{1}{2\pi} d\tau \mu(d\xi)$, i.e.
$$E[\cM(A)\overline{\cM(B)}]=\Pi(A \cap B) \quad \mbox{for any} \ A,B\in \cB_{b}(\bR^2).$$
Note that we have denoted by $\cB_b(\bR^2)$ the set of bounded Borel sets of $\bR^2$.

For any $t\geq 0$ and $\varphi \in \cD(\bR)$, we denote $X_{t}(\varphi)=X(1_{[0,t]}\varphi)$. This is well-defined, since $1_{[0,t]}\varphi \in \cH$ (see page 1128 of \cite{B12-Fourier}). Moreover,
$$X_t(\varphi)=\int_{\bR}\cF \varphi(\xi)M_t(d\xi),$$
where $\{M_t(A);A \in \cB_{b}(\bR)\}$ is a Gaussian complex random measure on $\bR$ with zero-mean and control measure $t \mu$, given by
\begin{equation}
M_t(A)=\int_{\bR^2}\frac{1-e^{-i\tau t}}{i\tau}1_{A}(\xi)\cM(d\tau,d\xi).
\label{eq:333}
\end{equation}
The process $M=\{M_t(A);t \geq 0,A \in \cB_b(\bR)\}$ is a {\em martingale measure} (as defined in \cite{walsh86}, but with complex values), with respect to the filtration
$$\cF_t=\sigma\{X_s(\varphi);s \in [0,t],\varphi \in \cD(\bR)\} \vee \cN, \quad t \geq 0$$
(see Appendix \ref{appendix-mart}). Here we denote by $\cN$ the class of $P$-negligible sets in $\cF$.
As we will show in the next section, the martingale measure $M$ plays an important role in the construction of the stochastic integral with respect to $X$.

For any $\phi \in L_{\bC}^2(\bR,\mu)$, we can define the integral $M_t(\phi):=\int_{\bR}\phi(\xi)M_t(d\xi)$ as an element in $L_{\bC}^2(\Omega)$, by approximation with simple functions. For any $s,t>0$ and $\phi,\psi \in L_{\bC}^2(\bR,\mu)$,
\begin{equation}
\label{cov-M}
E[M_t(\phi)\overline{M_s(\psi)}]=(t \wedge s) \int_{\bR}\phi(\xi) \overline{\psi(\xi)}\mu(d\xi).
\end{equation}
Hence, $X_t(\varphi)=M_t(\cF \varphi)$ for any $\varphi \in \cD(\bR)$. For any $s<t$ and $\varphi,\psi \in \cD(\bR)$,
$$E[X_t(\varphi)X_s(\psi)]=(t \wedge s)\int_{\bR}\cF \varphi(\xi)\overline{\cF \psi(\xi)}\mu(d\xi).$$

\begin{remark}
\label{random-field}
{\rm A random field $\{X(t,x);\, t\geq 0,x \in \bR\}$ can be
naturally associated to our noise $X$. In fact,
using an approximation argument it can be shown that for any $t \geq 0$ and $x\in \bR$, one has
$1_{(0,t] \times (0,x]} \in \cH$ and
$X(t,x):=X(1_{(0,t] \times (0,x]})=\int_{\bR}\cF 1_{(0,x]}(\xi)M_t(d\xi)$
(this follows from \cite[Thm 3.10]{B12-Fourier}).
The process $\{X(t,x);\, t\geq 0,x \in \bR\}$ has wide-sense stationary increments
(in the sense of \cite{BGP12}), with covariance
$$E[X(t,x)X(s,y)]=(t \wedge s)\int_{\bR}\cF 1_{[0,x]}(\xi) \overline{\cF 1_{[0,y]}(\xi)}\,\mu(d\xi)$$
(see Theorem 2.7 of \cite{BGP12}).
Finally, we observe that, for any fixed $t>0$, the process  $\{X_t(\varphi);\varphi \in \cD(\bR)\}$
can be identified with the distributional derivative $\partial_{x}$ of $\{X(t,x)\}_{x \in \bR}$,
which  is consistent with the developments in Section 2.2 of \cite{jolis10}.}
\end{remark}

\subsection{The stochastic integral}
\label{section-stoch-int}

In this section, we construct the stochastic integral with respect to $X$.
This construction is essentially the same as the one described in Section 4 of \cite{BGP12},
where the authors develop an integral with respect to the random field $\{X(t,x); \,t \geq 0,x \in \bR\}$
defined in Remark \ref{random-field}. However, we have chosen to focus more on the
stationary random distribution
$X=\{X(\varphi);\varphi \in \cD((0,\infty) \times \bR)\}$, so that our presentation is fully
consistent with the mathematical framework in the theory of SPDEs (see, e.g., \cite{dalang99,PZ07}).
As already mentioned, here we still consider that $H\in (0,1)$.

\smallskip

Note that for any interval $(x,y] \subset \bR$, we can define the random variable
\begin{equation}
\label{def-X-wrt-M}
X_t((x,y]):=X(1_{(0,t] \times (x,y]})=\int_{\bR}\cF 1_{(x,y]}(\xi)M_t(d\xi).
\end{equation}
However, in general, for an arbitrary set $A \in \cB_b(\bR)$, $X_t(A)$ can not be defined by the relation $X_t(A):=X(1_{(0,t] \times A})$, since $1_{(0,t] \times A}$ may not be in $\cH$. (Theorem 3.10 of \cite{B12-Fourier} can only be applied for sets $A$
for which $\int_{\bR} |\cF 1_{A}(\xi)|^2\mu(d\xi)<\infty$).
Precisely, this is the case whenever $H<\frac12$.
Therefore, in the general setting, the stochastic integral with respect to $X$ cannot be constructed using the approach of \cite{dalang99}.

The construction of the integral with respect to $X$ will be based on the random variables $X_t((x,y])$
given above. The properties of this integral are obtained indirectly,
using its relationship with the integral with respect to $M$. The integrals with respect to $M$
which will appear below are
defined as in Walsh's lecture notes \cite{walsh86}. We assume that the reader is familiar with
such kind of stochastic integrals.

\smallskip

In order to proceed with the construction of the integral with respect to $X$,
let us denote by $\cE_r$ the set of (real) linear combinations of processes of the form
\begin{equation}
\label{R-elem-process}
g(\omega,t,x)=Y(\omega)1_{(a,b]}(t) 1_{(u,v]}(x),
\end{equation}
where $0 \leq a<b$, $Y$ is a $\bR$-valued bounded $\cF_a$-measurable random variable, and $u,v \in \bR$ with $u<v$. The subscript $r$ in $\cE_r$ emphasizes the fact that a ``rectangle'' (of form $(u,v]$) is used in \eqref{R-elem-process}. If $g \in \cE_r$ is of form \eqref{R-elem-process}, we define the stochastic integral of $g$ with respect to $X$ by:
$$(g \cdot X)_t= Y(X_{t \wedge b}((u,v])-X_{t \wedge a}((u,v])).$$
This definition is extended by linearity to all processes in $\cE_r$.
The stochastic integral $(g \cdot X)_t$ can be expressed as an integral
with respect to the martingale measure $M$, as follows.

\begin{lemma}
\label{integral-elementary}
For any $g \in \cE_r$, the variable $(g \cdot X)_t$ has the spectral representation:
\begin{equation}
\label{spectral-rep-elem}
(g \cdot X)_t=\int_0^t \int_{\bR} \cF g(s,\cdot)(\xi)M(ds,d\xi).
\end{equation}
\end{lemma}

\noindent {\bf Proof:} It is enough to assume that $g$ is of form \eqref{R-elem-process}. The general case follows by linearity.
Using \eqref{def-X-wrt-M}, it follows that
$(g \cdot X)_t=  Y(M_{t \wedge b}(\psi)-M_{t \wedge a}(\psi))$, where $\psi:=\cF 1_{(u,v]} \in L_{\bC}^2(\bR,\mu)$.
On the other hand, the process
$$\Phi(\omega,t,x):=\cF g(\omega,t,\cdot)(\xi)=Y(\omega) 1_{(a,b]}(t) \psi(\xi)$$
is in $L_{\bC}^2(M)$, but is {\em not} an elementary process. By approximating $\psi$ with a sequence $(\psi_n)_n$ of simple functions, it can be shown that
$$\int_0^t \int_{\bR} \Phi (s,\xi)M(ds,d\xi)=Y(M_{t \wedge b}(\psi)-M_{t \wedge a}(\psi)).$$
The conclusion follows. \qed

\vspace{3mm}

Fix $T>0$. Similarly to \cite{dalang99}, we let $\cP_0$ be the completion of $\cE_r$ with respect to $\|\cdot\|_{0}$, where
$$\|g\|_{0}^2=E\int_0^T  \int_{\bR} |\cF g(t,\cdot)(\xi)|^2\mu(d\xi)dt.$$
In view of \eqref{spectral-rep-elem} and the isometry property of
Walsh's stochastic integral with respect to $M$, the map $g \mapsto \{(g \cdot X)_t\}_{t \in [0,T]}$ is an
isometry between $\cE_r$ and a subspace of the space of continuous square-integrable martingales $(M_t)_{t \in [0,T]}$ with $M_0=0$, endowed with the norm $\|M\|=\{E(M_T^2)\}^{1/2}$. This map can be extended to $\cP_0$. We denote the image of $g \in \cP_0$ under this map by
$$(g \cdot X)_t=\int_0^t \int_{\bR}g(s,x)X(ds,dx), \quad t\in [0,T].$$

We now identify the elements of the space $\cP_0$.

\begin{definition}
\label{pred-def-Sprime}
{\rm
We say that a function $S: \Omega \times [0,T] \to \cS'(\bR)$ is {\em predictable} if the map $(\omega,t) \mapsto S(\omega,t)(\varphi)$ is predictable, for any $\varphi \in \cS(\bR^d)$.
We will denote by $\cP_{\Omega \times \bR_+}$ the predictable $\sigma$-field on
$\Omega\times \bR_+$.
}
\end{definition}

\begin{remark}
\label{pred-func-Sprime}
{\rm
If $S: \Omega \times [0,T] \to \cS'(\bR)$ coincides with a function
$g:\Omega \times [0,T] \times \bR \to \bR$ (i.e. $S(\omega,t)(\varphi)=\int_{\bR^d}
g(\omega,t,x)\varphi(x)dx$ for all $\varphi \in \cS(\bR)$) and $g$ is predictable,
then $S$ is predictable (in the sense of Definition \ref{pred-def-Sprime}).
This follows by Fubini's theorem.
}
\end{remark}

\begin{remark}
\label{Fourier-meas}
{\rm
If $S: \Omega \times [0,T] \to \cS'(\bR^d)$ is a predictable function
such that $\cF S(\omega,t,\cdot)$ is a function for all $(\omega,t)$, then
by Lemma 4.2 of \cite{BGP12}, there exists a $\cP_{\Omega \times \bR_+} \times \cB(\bR)$-measurable function $\Phi: \Omega \times [0,T] \times \bR \to \bC$ such that for all $(\omega,t)$,
$$\cF S(\omega,t,\cdot)(\xi)=\Phi(\omega,t,\xi) \quad \mbox{for almost all} \ \xi \in \bR.$$
Below we will work with $\Phi(\omega,t,\xi)$, but we will write $\cF S(\omega,t,\cdot)(\xi)$.
}
\end{remark}

We consider the set $\Lambda_X$ of predictable functions $S:\Omega \times [0,T] \to \cS'(\bR)$ such that $\cF S(\omega,t,\cdot)$ is a locally integrable function for any $(\omega,t)$ and
$$E\int_0^T \int_{\bR} |\cF S(t,\cdot)(\xi)|^2 \mu(d\xi)dt<\infty.$$
The space $\Lambda_X$ is endowed with the inner product
$$\langle S_1, S_2 \rangle_{\Lambda_X}=E \int_0^T \int_{\bR} \cF S_1(t,\cdot)(\xi)\overline{S_2(t,\cdot)(\xi)}\mu(d\xi)dt.$$
We let $\|S\|_{\Lambda_X}^2=\langle S,S\rangle_{\Lambda_X}$. We identify $S_1$ and $S_2$ if $\|S_1-S_2\|_{\Lambda_X}=0$.

\begin{theorem}
\label{space-LambdaX}
The set $\Lambda_{X}$ coincides with $\cP_0$. For any $S \in \Lambda_X$ and $t \in [0,T]$,
\begin{equation}
\label{isometry}
E\left|\int_0^{t}\int_{\bR}S(s,x)X(ds,dx)\right|^2=E\int_0^{t}\int_{\bR} |\cF S(s,\cdot)(\xi)|^2\mu(d\xi)ds,
\end{equation}
and $(S \cdot X)_t$ admits the spectral representation:
\begin{equation}
\label{equality-intXM}
\int_0^{t}\int_{\bR}S(s,x)X(ds,dx)=\int_0^{t}\int_{\bR} \cF S(s,\cdot)(\xi)M(ds,d\xi) \quad {a.s.}
\end{equation}
For any $S \in \Lambda_X$, the predictable quadratic variation of $S \cdot X$ is:
\begin{equation}
\label{quadr-var-X}
\langle S \cdot X \rangle_t=\int_0^t \int_{\bR}|\cF S(s,\cdot)(\xi)|^2\mu(d\xi)ds, \quad t \in [0,T].
\end{equation}
\end{theorem}

\noindent {\bf Proof:} Clearly, $\langle g,h \rangle_{0}=\langle g,h\rangle_{\Lambda_X}$
for any $g,h \in \cE_r$, where we have denoted by
$\langle\cdot,\cdot\rangle_0$ the inner product associated with the norm
$\|\cdot\|_0$. By Theorem 4.3 of \cite{BGP12}, we know that $\cE_r$ is dense in $\Lambda_X$, and $\Lambda_X$ is complete. Hence, $\cP_0=\Lambda_X$.
Relation \eqref{equality-intXM} follows by an approximation argument, using Lemma \ref{integral-elementary} and the fact that $\cE_r$ is dense in $\Lambda_X$.
Relation \eqref{quadr-var-X} follows from \eqref{equality-intXM}. \qed

\subsection{A criterion for integrability}
\label{section-new-char}

In this section, we obtain a new criterion for integrability with respect to $X$,
which plays a crucial role in the present article. Here, we assume that $H\in (0,\frac12)$.

Throughout this article, we say that a measurable function $g:\bR \to \bR$ is {\em tempered} if
there exists a tempered distribution $T_g\in \cS'(\bR)$ such that
$T_g \varphi=\int_{\bR} g(x)\varphi(x)dx$, for all $\varphi\in \cS(\bR)$.

If $g$ is a tempered function, the Fourier transform of $g$ in $\cS'(\bR)$ is a tempered distribution,
defined by $\cF g (\phi)=\int_{\bR} g(x)\cF \phi(x)dx, \phi \in \cS(\bR)$. When this distribution is a locally integrable function, denoted also by $\cF g$, we have
$$\int_{\bR}g(x)\cF \phi(x)dx=\int_{\bR}\cF g(\xi)\phi(\xi)d\xi \quad \mbox{for all} \ \phi \in \cS(\bR).$$
In this case, the function $\cF g$ is also tempered.

We begin with a deterministic result, related to the theory of fractional Sobolev spaces,
which indeed slightly improves Proposition 3.4 in \cite{DPV}.

\begin{proposition}
\label{propA}
Let $g:\bR \to \bR$ be a tempered function whose Fourier transform in $\cS'(\bR)$ is a locally integrable function. For any $0<H<1/2$,
\begin{equation}
\label{key-identity}
c_H\int_{\bR}|\cF g(\xi)|^2 |\xi|^{1-2H}d\xi=C_H \int_{\bR^2} |g(x)-g(y)|^2|x-y|^{2H-2}dxdy,
\end{equation}
when either one of the two integrals above is finite. Here
$c_H$ is the constant given by \eqref{def-cH} and $C_H=H(1-2H)/2$.
\end{proposition}

\noindent {\bf Proof:}
First, assume that the integral on right-hand side of \eqref{key-identity}
is finite.
Then
\begin{align*}
\infty>\int_{\bR}\Big(\int_{\bR}\dfrac{|g(z+y)-g(y)|^2}{|z|^{2-2H}}\,dy\Big)dz & =
\int_{\bR}\,\left\|\dfrac{g(z+\cdot)-g(\cdot)}{|z|^{1-H}}\right\|_{L^2(\bR)}^2\,dz\\
&
=\frac1{2\pi}\int_{\bR}\,\left\|\cF\Big(\dfrac{g(z+\cdot)-g(\cdot)}{|z|^{1-H}}\Big)\right\|_{L^2(\bR)}^2\,dz,
\end{align*}
where we used Plancherel's theorem for the last equality. The application of Plancherel's theorem is
justified because the function $g(z+\cdot)-g(\cdot)$ belongs to $L^2(\bR)$,
for almost all $z \in \bR$. Since $\cF g$ is a tempered function, taking into account the definition and
properties of the Fourier transform of in $\cS'(\bR)$, we can infer that
$$\big(\cF g(z+\cdot)\big)(\xi)=e^{iz\xi}\,\cF g(\xi)\quad \xi\text{-a.e.},$$
and therefore, the last expression is equal to
\begin{align*}
&\frac1{2\pi}\int_{\bR}\int_{\bR} \dfrac{|e^{i\xi
z}-1|^2}{|z|^{2-2H}}\,|\cF g(\xi)|^2\,d\xi\,dz=\frac{\Gamma(2H+1)\sin(\pi H)}{\pi H(1-2H)}\int_{\bR}|\cF
g(\xi)|^2|\xi|^{1-2H}d\xi,
\end{align*}
using Fubini's theorem and Lemma \ref{elem-lem2} (see Appendix \ref{appendix-elem}). This proves that the integral on the left-hand side of \eqref{key-identity} is finite.
The constant which appears in front of the last integral above is exactly $c_H/C_H$.

When the integral on the left-hand side of \eqref{key-identity} is finite, we can use the same argument as above, but in reverse order. \qed

\vspace{3mm}

Based on the previous result, we now identify a subset of $\Lambda_{X}$.

\begin{theorem}
\label{theoremA}
Let $S:\Omega \times [0,T] \times \bR \to \bR$ be a predictable function,
such that
for almost all $(\omega,t) \in \Omega \times [0,T]$, $S(\omega,t,\cdot)$ is a tempered function whose Fourier transform $\cF S(\omega,t,\cdot)$ in $\cS'(\bR)$ is a locally integrable function. If
\begin{equation}
\label{IT-finite}
I(T):=C_HE\int_0^T \int_{\bR} \int_{\bR}|S(t,x)-S(t,y)|^2 |x-y|^{2H-2}dxdydt<\infty
\end{equation}
then $S \in \Lambda_X$ and $E|(S \cdot X)_T|^2=I(T)$. Moreover, for any $p \geq 2$,
\begin{equation}
\label{BDG}
E|(S \cdot X)_T|^p \leq z_p C_{H}^{p/2} E \left(\int_0^T \int_{\bR} \int_{\bR} |S(t,x)-S(t,y)|^2 |x-y|^{2H-2} \,dx \,dy \,dt \right)^{p/2}
\end{equation}
where $z_p$ is the constant in the Burkholder-Davis-Gundy inequality for continuous martingales.
\end{theorem}

\noindent {\bf Proof:} Note that the function
$(\omega,t,x,y) \mapsto |S(\omega,t,x)-S(\omega,t,y)|^2 |x-y|^{2H-2}$ is $\cP_{\Omega \times \bR_{+}}
\times \cB(\bR^2)$-measurable, due to the fact that
\[
 \cP_{\Omega\times \bR_+\times \bR} \subset \cP_{\Omega\times \bR_+}  \times \cB(\bR).
\]
By \eqref{IT-finite}, there exists a set $N \in \cP_{\Omega \times \bR_{+}}$ with $(P\times {\rm Leb})(N)=0$ such that for all $(\omega,t) \not\in N$,
$$\int_{\bR} \int_{\bR}|S(\omega,t,x)-S(\omega,t,y)|^2 |x-y|^{2H-2}dxdy<\infty.$$

We apply Proposition \ref{propA} to the function $g=S(\omega,t,\cdot)$ with $(\omega,t) \not\in N$. We obtain that for any $(\omega,t) \not \in N$.
\begin{equation}
\label{eq1}
c_H \int_{\bR}|\cF S(\omega,t,\cdot)(\xi)|^2 |\xi|^{1-2H}d\xi=C_H \int_{\bR}\int_{\bR}|S(\omega,t,x)-S(\omega,t,y)|^2|x-y|^{2H-2}dxdy.
\end{equation}

For any $(\omega,t) \in \Omega \times [0,T]$, we identify the function $S(\omega,t,\cdot)$ with the distribution in $\cS'(\bR)$ induced by it. By Remark \ref{pred-func-Sprime},  $S:\Omega \times [0,T] \to \cS'(\bR)$ is predictable (in the sense of Definition \ref{pred-def-Sprime}).

As in Remark \ref{Fourier-meas}, there exists a function $\Phi:\Omega \times [0,T] \times \bR \to \bC$ which is $\cP_{\Omega \times \bR_{+}} \times \cB(\bR)$-measurable, such that for any $(\omega,t)$,
$\cF S(\omega,t,\cdot)(\xi)=\Phi(\omega,t,\xi)$ for almost all $\xi \in \bR$.
Hence for any $(\omega,t) \in \Omega \times [0,T]$,
\begin{equation}
\label{eq2}
\int_{\bR}|\cF S(\omega,t,\cdot)(\xi)|^2 |\xi|^{1-2H}d\xi=\int_{\bR}|\Phi(\omega,t,\xi)|^2 |\xi|^{1-2H}d\xi.
\end{equation}
From \eqref{eq1} and \eqref{eq2}, it follows that for any $(\omega,t) \not\in N$,
\begin{equation}
\label{Phi-equality}
c_H\int_{\bR}|\Phi(\omega,t,\xi)|^2 |\xi|^{1-2H}d\xi=C_H \int_{\bR}\int_{\bR}|S(\omega,t,x)-S(\omega,t,y)|^2|x-y|^{2H-2}dxdy.
\end{equation}

We now take the integral with respect to $P(d\omega)dt$. We obtain that
\begin{equation}
\label{It-equality}
c_HE\int_{0}^T \int_{\bR}|\cF S(t,\xi)|^2 \,|\xi|^{1-2H}d\xi dt=I(T)<\infty.
\end{equation}
This proves that $S \in \Lambda_X$. The fact that $E|(S \cdot X)_T|^2=I(T)$ follows from \eqref{isometry} and \eqref{It-equality}.

We now prove \eqref{BDG}. By Burkholder-Davis-Gundy inequality,
\begin{equation}
\label{BDG2}
E|(S \cdot X)_T|^p \leq z_p E \left(\langle S \cdot X \rangle_{T}^{p/2}\right).
\end{equation}
By \eqref{quadr-var-X}, we know that for almost all $\omega\in \Omega$,
\begin{equation}
\label{quadr-var1}
\langle S \cdot X \rangle_{T}(\omega) = c_H \int_{0}^{T} \int_{\bR}|\Phi(\omega,t,\xi)|^2 \,|\xi|^{1-2H}d\xi dt.
\end{equation}

Let $F$ be the set of $(\omega,t)$'s for which \eqref{Phi-equality} does not holds, and $F_{\omega}=\{t \in [0,T];(\omega,t) \in F\}$. Since $(P \times {\rm Leb})(F)=0$, by Fubini's theorem, ${\rm Leb}(F_{\omega})=0$ for almost all $\omega$. Hence, there exists a set $\Omega_0$ with $P(\Omega_0)=1$ such that for all $\omega \in \Omega_0$ fixed, equality \eqref{Phi-equality} holds for almost all $t \in [0,T]$. Taking the integral with respect to $dt$, we obtain that for any $\omega \in \Omega_0$,
\begin{align}
\nonumber
&c_H \int_{0}^{T} \int_{\bR}|\Phi(\omega,t,\xi)|^2 \,|\xi|^{1-2H}d\xi dt \\
\label{quadr-var2}
& \qquad \quad = C_H \int_0^T \int_{\bR} \int_{\bR} |S(\omega,t,x)-S(\omega,t,y)|^2 |x-y|^{2H-2}\,dx \,dy \,dt.
\end{align}
From \eqref{quadr-var1} and \eqref{quadr-var2}, we infer that:
$$\langle S \cdot X \rangle_{T} = C_H \int_0^T \int_{\bR} \int_{\bR} |S(t,x)-S(t,y)|^2 |x-y|^{2H-2}\,dx \,dy \,dt \quad {\rm a.s.}$$
Relation \eqref{BDG} follows, using\eqref{BDG2}. \qed

\section{Proof of Theorem \ref{main-thm}}
\label{section-proof}

This section is dedicated to the proof of Theorem \ref{main-thm}. In particular,
from now on, we assume that $H\in(\frac14,\frac12)$.
In Section \ref{section-prelim}, we gather some preliminary (deterministic) results which
are needed in the subsequent sections.
Section \ref{section-Picard} will be devoted to prove that the sequence $(u^n)_{n \geq 0}$
of Picard iterations is well-defined, for both equations \eqref{wave} and \eqref{heat}.
In Section \ref{section-conv-Picard}, we show that the sequence $(u^n)_{n \geq 0}$ converges
(in a certain Banach space, which is defined below), and its limit is the desired solution.
Moreover, we show that the solution is unique. Finally, in Section \ref{section-gronwall},
we state and prove an extension of Gronwall's lemma, which is of independent interest,
and is used in Section \ref{section-conv-Picard}.

\subsection{Preliminary results}
\label{section-prelim}

In this section, we give some estimates for various integrals containing the fundamental solution $G$ of the wave or heat equations.
Recall that the Fourier transform of $G_t$ is:
$$\cF G_t(\xi)=\frac{\sin(t|\xi|)}{|\xi|}\quad \mbox{for the wave equation},$$
$$\cF G_t(\xi)=\exp\left(-\frac{t|\xi|^2}{2}\right) \quad \mbox{for the heat equation}.$$

\begin{lemma}
\label{lemma-G1} Let $G$ be the fundamental solution of the wave or
heat equation. Then, for both equations, the integral
$$A_T(\alpha):=\int_0^T\int_{\bR}|\cF G_t(\xi)|^2\,|\xi|^{\alpha}\,d\xi\,dt$$
converges if and only if $\alpha \in (-1,1)$. When the integral
converges, we have
\begin{eqnarray}
A_T(\alpha)&=&  2^{1-\alpha}\,C_{\alpha }\frac{1}{2-\alpha} \,T^{2-\alpha}
\quad \mbox{for wave equation}, \label{eq:77}\\
A_T(\alpha)&=&\frac{2}{1-\alpha}\,\Gamma\left(\frac{\alpha+1}{2}\right)
T^{(1-\alpha)/2} \quad \mbox{for heat equation}, \label{eq:78}
\end{eqnarray}
where
$$C_{\alpha}=
\left\{
\begin{array}{ll}
(1-\alpha)^{-1}\Gamma(\alpha) \sin(\pi \alpha/2) & \mbox{if $\alpha \in (0,1)$} \\
\alpha^{-1}(1-\alpha)^{-1}\Gamma(1+\alpha)\sin(\pi \alpha/2) & \mbox{if $\alpha \in (-1,0)$}\\
\pi/2 & \mbox{if $\alpha=0$}.
\end{array} \right.
$$
\end{lemma}

\noindent {\bf Proof:} For the wave equation, we use Lemma
\ref{elem-lem1} with $\alpha'=1-\alpha \in (0,2)$:
\begin{align}
\nonumber
&\int_{\bR} |\cF G_t(\xi)|^2|\,\xi|^{\alpha}d\xi=\int_{\bR}\frac{\sin^2(t|\xi|)}{|\xi|^{2-\alpha}}
d\xi=
2t^{1-\alpha}\int_0^{\infty}\frac{\sin^2 x}{x^{2-\alpha}}dx\\
\label{int-wave} & \quad
=t^{1-\alpha}\int_{0}^{\infty}\frac{1-\cos(2x)}{x^{2-\alpha}}dx= t^{1-\alpha}\, 2^{1-\alpha}C_\alpha,
\end{align}
where the last integral converges if and only if $\alpha \in
(-1,1)$. Thus, we obtain \eqref{eq:77}. For the heat equation,
\begin{eqnarray}
\nonumber
\int_{\bR} |\cF G_t(\xi)|^2\,|\xi|^{\alpha}d\xi &=& \int_{\bR} e^{-t|\xi|^2}|\xi|^{\alpha}d\xi=2\int_0^{\infty}e^{-t x^2} x^{\alpha}dx\\
\label{int-heat} & = &
t^{-(\alpha+1)/2}\Gamma\left(\frac{\alpha+1}{2}\right),
\end{eqnarray}
using the change of variable $y=tx^2$. The last integral converges
if and only if $\alpha>-1$. Finally, the integral $\int_0^T
t^{-(\alpha+1)/2}dt$ converges if and only if $\alpha<1$, whence we
also deduce \eqref{eq:78}. \qed

\medskip

\begin{remark}
\label{remark-H}
{\rm In the sequel, we will apply Lemma \ref{lemma-G1} with $\alpha=2(1-2H)$ (which imposes the restriction $H>1/4$), and also with $\alpha=1-2H$.}
\end{remark}

\begin{lemma}
\label{lemma-G2}
Assume that $1/4<H<1/2$. Let $G$ be the fundamental solution of the wave or heat equation. For any $0\leq a\leq b\leq T$, set
$$F(a,b):= \int_a^b \!\!\int_{\bR} G_{b-s}^2(z) \int_{\bR}
  |\cF G_{s-a} (\xi)|^2 |\xi|^{2(1-2H)}d\xi dz ds.$$
Then
$$F(a,b)=\left\{
\begin{array}{ll}
C(b-a)^{4H+1} & \mbox{for the wave equation}\\
C(b-a)^{2H-1} &\mbox{for the heat equation}
\end{array}
\right.$$
where $C>0$ is a constant depending on $H$, which is different for the two equations.
\end{lemma}

\noindent {\bf Proof}: Using \eqref{int-wave} and \eqref{int-heat} with $\alpha=2(1-2H)$, we see that
$$\int_{\bR} |\cF G_{s-a} (\xi)|^2 |\xi|^{2(1-2H)}d\xi =
\left\{
\begin{array}{ll}
C (s-a)^{4H-1} & \mbox{for the wave equation}\\
C(s-a)^{2H-3/2} & \mbox{for the heat equation}
\end{array}
\right.$$
where $C>0$ is a constant depending on $H$, which is different for the two equations.
This integral is finite because $H>1/4$. On the other hand,
\begin{equation}
\label{L2normG}
\int_{\bR} G_{t}(z)^2 dz=
\left\{
\begin{array}{ll}
2^{-1}t & \mbox{for the wave equation}\\
2\pi^{1/2}t^{-1/2} & \mbox{for the heat equation}
\end{array}
\right.
\end{equation}
Therefore, for the wave equation, we have:
\begin{eqnarray*}
F(a,b)&=&C 2^{-1}\int_a^b (b-s) (s-a)^{4H-1}ds =C2^{-1}\int_{0}^{b-a} (b-a-r)\, r^{4H-1}\, dr\\
&=& C 2^{-1}\beta (2,4H) (b-a)^{4H+1},
\end{eqnarray*}
where $\beta(a,b)=\Gamma(a)\Gamma(b)/\Gamma(a+b)$ is the beta function.
For the heat equation,
\begin{eqnarray*}
F(a,b)&=&C2\pi^{1/2}\int_a^b (b-s)^{-1/2} (s-a)^{2H-3/2}ds \\
&=& C2\pi^{1/2}\beta\left(\frac{1}{2}\, ,2H-\frac{1}{2}\right)(b-a)^{2H-1},
\end{eqnarray*}
using again the fact that $H>1/4$. \qed

\vspace{3mm}

\begin{lemma}
\label{lemma-G3}
Let $G$ be the fundamental solution of the wave or heat equation. For any $\alpha \in (-1,1)$ and for any $h \in \bR$,
$$\int_0^T \int_{\bR} (1-\cos(\xi h)) \, |\cF G_t(\xi)|^2 \,|\xi|^{\alpha}\, d\xi \,dt \leq
\left\{
\begin{array}{ll}
C T |h|^{1-\alpha} & \mbox{for the wave equation} \\
C |h|^{1-\alpha} & \mbox{for the heat equation},
\end{array}
\right.
$$
where $C=\int_{\bR} (1-\cos \eta)\eta^{\alpha-2}d\eta$.
\end{lemma}

\noindent {\bf Proof:} It is enough to consider the case $h>0$. If $h<0$, we use the fact that $\cos(\xi h)=\cos(\xi|h|)$. For the wave equation, since $\sin^2(t|\xi|) \leq 1$,
$$\int_0^T \int_{\bR} (1-\cos(\xi h)) \frac{\sin^2(t|\xi|)}{|\xi|^2}|\xi|^{\alpha} \,d\xi \,dt \leq T \int_{\bR}\frac{1-\cos(\xi h)}{|\xi|^{2-\alpha}}\,d\xi=C T h^{1-\alpha},$$
using Lemma \ref{elem-lem1} with $\alpha'=1-\alpha$. For the heat
equation,
$$\int_0^T \int_{\bR} (1-\cos(\xi h)) \, e^{-t|\xi|^2} |\xi|^{\alpha} \, d\xi \,dt =\int_{\bR}
(1-\cos (\xi h)) \frac{1-e^{-T|\xi|^2}}{|\xi|^2}|\xi|^{\alpha}d\xi
\leq C h^{1-\alpha},$$ using again Lemma \ref{elem-lem1} and the
fact that $1-e^{-T|\xi|^2} \leq 1$. \qed

\begin{lemma}
\label{lemma-G4}
Let $G$ be the fundamental solution of the wave or heat equation. For any $\alpha \in (-1,1)$ and for any $h\in \bR$
$$\int_{0}^{T}\int_{\bR}|\cF G_{t+h}(y)-\cF G_t(y)|^2 \,|\xi|^{\alpha}\,d\xi \,dt\leq
\left\{
\begin{array}{ll}
C T |h|^{1-\alpha} & \mbox{for the wave equation} \\
C |h|^{(1-\alpha)/2} & \mbox{for the heat equation},
\end{array}
\right.$$
where $C$ is a constant depending on $\alpha$ (which is different for the two equations).
\end{lemma}

\noindent {\bf Proof:} We suppose that $h>0$. The case $h<0$ is similar.
For the wave equation, using the fact that
$|\sin((t+h)|\xi|)-\sin(t|\xi|)|^2 \leq C \min(1,h|\xi|)^2$,
we see that the integral is smaller than
$$C\, T \int_{\bR} \frac{\min(1,h|\xi|)^2}{|\xi|^{2-\alpha}}\,d\xi= C\,
T \, h^{1-\alpha} \int_{\bR}\frac{\min(1,|\eta|)^2}{|\eta|^{2-\alpha}}\,d\eta \leq CT \, h^{1-\alpha},$$
using the change of variables $\eta=h\xi$.
For the heat equation, the integral is
\begin{eqnarray*}
& & \int_0^t\int_{\bR}e^{-s|\xi|^2} (1-e^{-h|\xi|^2/2})^2\,|\xi|^{\alpha}\,d\xi\,ds
=\int_{\bR}(1-e^{-h|\xi|^2/2})^2 \, \frac{1-e^{-t|\xi|^2}}{|\xi|^2} \,|\xi|^{\alpha}\,d\xi \\
 & & \quad \leq
\int_{\bR} \, \frac{(1-e^{-h|\xi|^2/2})^2}{|\xi|^{2-\alpha}}\,d\xi=
h^{(1-\alpha)/2}\int_{\bR}\frac{(1-e^{-\eta^2/2})^2}{|\eta|^{2-\alpha}}d\eta,
\end{eqnarray*}
using the change of variables $\eta=h^{1/2}\xi$. The last integral
is seen to be finite using the fact that $1-e^{-\eta^2/2} \leq
\eta^2/2$ for $|\eta| \leq 1$. \qed

\subsection{Picard iteration scheme}
\label{section-Picard}

In this section, we show that the sequence $(u^n)_{n \geq 0}$ of the
Picard iterations is well-defined and satisfies some properties.
This sequence is defined iteratively. We set $u^0(t,x)=w(t,x)$ for
any $t \in [0,T]$ and $x \in \bR$. For any $n \geq 0$ and for any
$(t,x) \in [0,T] \times \bR$, we let
\begin{equation}
\label{def-Picard}
u^{n+1}(t,x)=w(t,x)+\int_0^{t}\int_{\bR} G_{t-s}(x-y)\,\sigma(u^{n}(s,y))\,X(ds,dy),
\end{equation}
provided that the stochastic integral is well-defined, in the sense
explained in Section \ref{section-stoch-int} above. For the next
result, the function $\sigma$ does not need to be affine.
Recall that the term $w(t,x)$ corresponds to the contribution of the initial data, and
explicit formulas for it are given in the Introduction.

\begin{theorem}
\label{Picard-well-def} Let $\sigma$ be an arbitrary Lipschitz
function. Let $p \geq 2$ be fixed. Then, for any $n \geq 0$,
\begin{equation}
\left.\begin{array}{rcl}
& & \displaystyle u^n(t,x) \ \mbox{is well-defined for any} \ (t,x) \in [0,T] \times \bR,  \\[2ex]
& & \displaystyle \sup_{(t,x) \in [0,T] \times \bR}E|u^n(t,x)|^p<\infty, \quad \mbox{and} \\[1ex]
& & \displaystyle \sup_{(t,x) \in [0,T] \times \bR} \int_0^t \int_{\bR^2}G_{t-s}^2(x-y)\frac{\Big(E|u^n(s,y)-u^n(s,z)|^p\Big)^{2/p}}{|y-z|^{2-2H}}\,dy \,dz \,ds<\infty
\end{array}\right\} \label{propertyP} \tag{P}
\end{equation}
and, for any $h\in \bR$ with $|h|<1$,
\begin{equation}
\left.\begin{array}{rcl}
& & \displaystyle \sup_{(t,x) \in [0,T] \times \bR}|u^n(t,x+h)-u^n(t,x)|^2 \leq C_n |h|^{2H} \\
& & \displaystyle \sup_{(t,x) \in [0,T \wedge (T-h)] \times
\bR}|u^n(t+h,x)-u^n(t,x)|^2 \leq C_n |h|^{\beta},
\end{array} \right\} \label{propertyQ} \tag{Q}
\end{equation}
where $\beta=2H$ for the wave equation, and $\beta=H$ for the heat
equation. Here $C_n$ is a constant which depends on $n$ (and also on
$H,T,\sigma,u_0$ and $v_0$).
\end{theorem}

\noindent {\bf Proof:} By induction, we prove that properties
\eqref{propertyP} and \eqref{propertyQ} hold for any $n \geq 0$. As
already mentioned, the constant $C$ (depending on $p,H,T,\sigma,u_0$
and $v_0$) may be different from line to line. We split the proof in
four steps, as follows.

\smallskip

\noindent {\it Step 1.} We start by checking properties
\eqref{propertyP} and \eqref{propertyQ} for $n=0$. It is clear that
the variable $u^0(t,x)=w(t,x)$ is well-defined for any $(t,x)$.
Using the particular form of $w$ for each equation, and the fact
that $u_0$ and $v_0$ are bounded, we see immediately that for both
equations,
\begin{equation}
\label{w-bounded}
\sup_{(t,x) \in [0,T] \times \bR}|w(t,x)|<\infty.
\end{equation}
Let us postpone for the moment the proof of the third condition in
\eqref{propertyP}, so that we first check the validity of property
\eqref{propertyQ} for the case $n=0$.

For this, we consider separately the wave and heat equations, using
the corresponding formulas for $w(t,x)$ for each equation. For the
wave equation, $|w(t,x+h)-w(t,x)|^2 \leq 2(A_1+A_2)$, where
\begin{align*}
& A_1 =\left|\int_{\bR}(G_t(x+h-y)-G_t(x-y))v_0(y)dy \right|^2\\
& A_2 =\frac{1}{4}\,|u_0(x+h+t)-u_0(x+t)+u_0(x+h-t)-u_0(x-t)|^2.
\end{align*}
Since $u_0$ is $H$-H\"older continuous, $A_2 \leq C|h|^{2H}$. By a
change of variables, and H\"older's inequality (using the fact that
$\int_{\bR}G(t,y)dy=t$), we have:
\begin{eqnarray*}
A_1 &=& \left|\int_{\bR}G_t(y)(v_0(x+h-y)-v_0(x-y)) dy \right|^2 \\
&\leq & t \int_{\bR}G_t(y)\,|v_0(x+h-y)-v_0(x-y)|^2 dy
 \leq  C t^2 |h|^{2H},
\end{eqnarray*}
where for the last inequality we used the fact that $v_0$ is
$H$-H\"older continuous. We now examine the time increments of $w$.
We suppose that $h>0$. The case $h<0$ is similar. We have
 $|w(t+h,x)-w(t,x)|^2 \leq 2(B_1+B_2)$, where
\begin{align*}
& B_1 =\left|\int_{\bR}(G_{t+h}(x-y)-G_t(x-y))v_0(y)dy \right|^2\\
& B_2 =\frac{1}{4}\,|u_0(x+t+h)-u_0(x+t)+u_0(x-t-h)-u_0(x-t)|^2.
\end{align*}
Since $u_0$ is $H$-H\"older continuous, $B_2 \leq C|h|^{2H}$. By
H\"older's inequality,
\begin{align*}
B_1&= \frac{1}{4}\left|\int_{\bR}(1_{\{x+t<y<x+t+h\}}+1_{\{x-t-h<y<x-t\}})\, v_0(y)\, dy\right|^2\\
&\leq \frac{1}{2}\,h
\,\int_{\bR}(1_{\{x+t<y<x+t+h\}}+1_{\{x-t-h<y<x-t\}})\, |v_0(y)|^2
\,dy \leq C h^2,
\end{align*}
where for the last inequality we used the fact that $v_0$ is
bounded.

In the case of the heat equation, the space increments of $w$ are
treated similarly to the term $A_1$ above (with $v_0$ replaced by
$u_0$). For the time increments, we assume again that $h>0$, the
case $h<0$ being similar. Using the semigroup property of $G$,
H\"older's inequality (since $\int_{\bR}G_t(y)dy=1$) and the fact
that $u_0$ is $H$-H\"older continuous, we obtain (see (4.8) of
\cite{sanz-sarra00}):
\begin{align*}
|w(t+h,x)-w(t,x)|^2 &=\left|\int_{\bR}G_{h}(y) \int_{\bR}G_{t}(x-z) (u_0(z-y)-u_0(z))dz dy\right|^2 \\
&\leq \int_{\bR}G_{h}(y) \int_{\bR}G_{t}(x-z) \, |u_0(z-y)-u_0(z)|^2 \,dz \, dy \\
& \leq C \int_{\bR}G_{h}(y)\, |y|^{2H}dy=C h^{H}.
\end{align*}
Therefore, we conclude that property \eqref{propertyQ} holds for
$n=0$.

It remains to show that the third condition in \eqref{propertyP}
holds in this case. Precisely, we write
$$\int_0^t \int_{\bR} G_{t-s}^2(x-y)\left( \int_{\bR}\frac{|w(s,y+z)-w(s,y)|^2}{|z|^{2-2H}}\,dz \right)\,dy \,ds=I'+I'',$$
where the terms $I'$ and $I''$ are obtained by splitting the $dz$
integral into two integrals, corresponding to the regions
$\{|z|>1\}$ and $\{|z| \leq 1\}$, respectively. Note that
$$I' \leq 4\, \sup_{(t,x)}|w(t,x)|^2 \int_0^t \int_{\bR}G_{t-s}^2(x-y) \left(\int_{|z|>1}\frac{1}{|z|^{2-2H}}\,dz\right) \,dy \,ds,$$
which is uniformly bounded for all $(t,x) \in [0,T] \times \bR$, due
to \eqref{L2normG}  and the fact that $H<1/2$. For $I''$, we use the
fact that $|w(s,y+z)-w(s,y)|^2 \leq C|z|^{2H}$, since we have
already proved that property \eqref{propertyQ} holds for $n=0$.
Hence,
$$I'' \leq \int_0^t \int_{\bR} G_{t-s}^2(x-y)\left( \int_{|z| \leq 1}|z|^{4H-2}\,dz \right)\,dy \,ds,$$
which is uniformly bounded in $(t,x)$, due to \eqref{L2normG} and the fact that $H>1/4$.
This concludes the first step of the proof.

\smallskip

From now on, we assume the following
induction hypothesis: properties \eqref{propertyP} and
\eqref{propertyQ} hold for $u^n$. We aim to prove that
\eqref{propertyP} and \eqref{propertyQ} also hold for $u^{n+1}$.

\smallskip

\noindent {\it Step 2.}
This step of the proof is devoted to show that $u^{n+1}(t,x)$ is
well-defined for all $(t,x)\in [0,T]\times \bR$ and it holds
\begin{equation}
\sup_{(t,x) \in [0,T] \times \bR}E|u^{n+1}(t,x)|^p<\infty.
\label{eq:79}
\end{equation}

In order to prove that  $u^{n+1}(t,x)$ is well-defined for any
$(t,x)\in [0,T] \times \bR$, we have to show that the stochastic
integral on the right-hand side of \eqref{def-Picard} is
well-defined. For this, we apply Theorem \ref{theoremA} to the
function
$$S_n(s,y)=G_{t-s}(x-y)\sigma(u^n(s,y))1_{[0,t]}(s)$$
(which depends also on $(t,x)$). We will show that:\\
{\em (i)} $u^n$ has a predictable modification (called also $u^n$);\\
{\em (ii)} $S_n(\omega,s,\cdot)$ is in $L^1(\bR)$ for almost all $(\omega,s) \in \Omega \times [0,T]$; \\
{\em (iii)} $S_n$ satisfies the following condition:
\begin{equation}
\label{finite-int-S}
\sup_{(t,x) \in [0,T] \times \bR}E\int_0^t \int_{\bR}\int_{\bR}|S_n(s,y)-S_n(s,z)|^2|y-z|^{2H-2}\,dy \,dz \,ds<\infty.
\end{equation}

To prove {\em (i)}, we note that $u^n$ is $L^2(\Omega)$-continuous
(by property \eqref{propertyQ}) and $u^n(t,x)$ is $\cF_t$-measurable for any $x \in \bR$ (by the construction of the stochastic integral). Therefore, $u^n$ has a predictable modification. (For this, we use the extensions to random fields of Theorem 30, Chapter IV of \cite{DM75} and Proposition 3.21 of \cite{PZ07}.)
To prove {\em (ii)}, we note that
$$E \int_0^T \int_{\bR} |S_n(s,y)|dyds \leq C \left(1+\sup_{(t,x)}E|u^n(t,x)|\right)
\int_0^T \int_{\bR} G_{t-s}(x-y)dyds<
\infty,$$
and hence, $\int_{\bR} |S_n(\omega,s,y)|dy<\infty$ for almost all $(\omega,s)$.
This proves {\em (ii)}. Note also that {\em (ii)} implies that the Fourier transform of
$S_n(s,\cdot)$, for a.a. $(\omega,s)$, is given by a function.

To prove {\em (iii)}, we bound the integral in \eqref{finite-int-S} by $2(I_1+I_2)$, where
\begin{eqnarray*}
I_1&=& E
\int_0^{t}\int_{\bR}\int_{\bR}G_{t-s}^2(x-y)
\dfrac{|\sigma(u^n(s,y))-\sigma(u^n(s,z))|^2}{|y-z|^{2-2H}}\,dy\,dz\,ds,\\
I_2&=& E
\int_0^{t}\int_{\bR} \int_{\bR} |\sigma(u^n(s,z))|^2
\dfrac{|G_{t-s}(x-y)-G_{t-s}(x-z)|^2}{|y-z|^{2-2H}}\,dy\,dz\,ds.
\end{eqnarray*}
Using the fact that $\sigma$ is Lipschitz, and Jensen's inequality,
we see that
$$I_1 \leq C\int_0^{t}\int_{\bR}\int_{\bR}G_{t-s}^2(x-y)
\dfrac{\Big(E|u^n(s,y)-u^n(s,z)|^p\Big)^{2/p}}{|y-z|^{2-2H}}\,dy\,dz\,ds,$$
which is uniformly bounded in $(t,x)$ by the induction hypothesis (the third condition in \eqref{propertyP}). By Jensen's inequality and Proposition \ref{propA},
\begin{eqnarray*}
I_2 & \leq & C\left(1+\sup_{(t,x)}E|u^n(t,x)|^2\right)
\int_0^{t}\int_{\bR}\int_{\bR}
\dfrac{|G_{t-s}(x-y)-G_{t-s}(x-z)|^2}{|y-z|^{2-2H}}\,dy\,dz\,ds \\
& \leq &  C\left[1+\left(\sup_{(t,x)}E|u^n(t,x)|^p\right)^{2/p}\right] \int_{0}^{t} \int_{\bR} |\cF G_{t-s}(\xi)|^{2}\,|\xi|^{1-2H}d\xi\, ds,
\end{eqnarray*}
the $\sup$ being finite by the induction hypothesis (the second property in \eqref{propertyP}). The last integral is uniformly bounded for $t \in [0,T]$, by Lemma \ref{lemma-G1}. This shows that $u^{n+1}(t,x)$ is well-defined.

Next, we show that \eqref{eq:79} holds. By the definition
\eqref{def-Picard} of $u^{n+1}(t,x)$ and Theorem \ref{theoremA},
$$E|u^{n+1}(t,x)|^p \leq C\left\{|w(t,x)|^p+E\left( \int_0^t \int_{\bR}\int_{\bR}\dfrac{|S_n(s,y)-S_n(s,z)|^2}{|y-z|^{2-2H}}\,dy \,dz \,ds  \right)^{p/2} \right\}.$$
The first term is uniformly bounded in $(t,x)$ by \eqref{w-bounded}. The second term is bounded by $C(J_1+J_2)$, where
\begin{eqnarray*}
J_1 & = & E \left(\int_0^t \int_{\bR^2}G_{t-s}^2(x-y)\dfrac{|\sigma(u^n(s,y))-\sigma(u^n(s,z))|^2} {|y-z|^{2-2H}}\,dy\,dz\, ds \right)^{p/2}
 \\
J_2  &= & E \left(\int_0^t\int_{\bR^2} |\sigma(u^n(s,z))|^2
\dfrac{|G_{t-s}(x-y)-G_{t-s}(x-z)|^2}{|y-z|^{2-2H}} \,dy\,dz\,ds\right)^{p/2}.
 \end{eqnarray*}
Using the fact that $\sigma$ is Lipschitz and applying Minkowski
inequality for integrals (see A.1, page 271 of \cite{stein70}), we see that
$$J_1 \leq C \left(\int_0^t \int_{\bR^2}G_{t-s}^2(x-y)\dfrac{\Big(E|u^n(s,y)-u^n(s,z)|^p\Big)^{2/p}}{|y-z|^{2-2H}} \,dy\,dz \,ds\right)^{p/2},$$
which is uniformly bounded in $(t,x)$, by the induction hypothesis
(the third condition in \eqref{propertyP}). Similarly for $J_2$, we have:
\begin{align*}
J_2 & \leq C \left(\int_0^t\int_{\bR^2} \Big(1+E|u^n(s,z)|^p\Big)^{2/p} \dfrac{|G_{t-s}(x-y)-G_{t-s}(x-z)|^2}{|y-z|^{2-2H}} \,dy\,dz\,ds\right)^{p/2}\\
& \leq C \left[1+\left(\sup_{(t,x)}E|u^n(t,x)|^p\right)^{2/p}\right]
\left( \int_{0}^{t}\int_{\bR^2}\frac{|G_{t-s}(x-y)-G_{t-s}(x-z)|^2}{|y-z|^{2-2H}}\,dy\,dz\,ds\right)^{p/2},
\end{align*}
the sup being bounded by the induction hypothesis.
As mentioned above, the last integral is uniformly bounded in $t \in [0,T]$.
Thus, we have obtained the validity of \eqref{eq:79}).

\smallskip

\noindent {\it Step 3.} Now, we prove that $u^{n+1}$ satisfies the
third condition in \eqref{propertyP}, i.e.
\begin{equation}
\label{bound2-n+1}
\sup_{(t,x)\in[0,T]\times\bR}\int_0^{t}\int_{\bR^2}
\frac{G_{t-s}(x-y)^2}{|z|^{2-2H}}
\Big(E|u^{n+1}(s,y+z)-u^{n+1}(s,y)|^p\Big)^{2/p}\,dy\,dz\,ds<\infty.
\end{equation}
By \eqref{def-Picard} and Theorem \ref{theoremA},
\begin{align*}
& E|u^{n+1}(s,y+z)-u^{n+1}(s,y)|^p \\
&\quad  \leq C |w(s,y+z)-w(s,y)|^p\\
&  \quad \quad + CE\left|\int_0^s \int_{\bR}(G_{s-r}(y+z-v)-G_{s-r}(y-v))\sigma(u^n(r,v))X(dr,dv)\right|^p\\
& \quad \leq C |w(s,y+z)-w(s,y)|^p\\
&  \quad \quad + C E\left(\int_0^s \int_{\bR^2}
|(G_{s-r}(y+z-v)-G_{s-r}(y-v))\sigma(u^n(r,v)) \right.\\
&  \quad \quad - \left.(G_{s-r}(y+z-\bar v)-G_{s-r}(y-\bar v))\sigma(u^n(r,\bar v))|^2\frac{1}{|v-\bar v|^{2-2H}}\,dv \,d\bar v \, dr\right)^{p/2}.
\end{align*}
Since we have already proved that the third condition in \eqref{propertyP} holds for $w$, it suffices to show that
\begin{align*}
&\int_0^{t}\int_{\bR^2}\frac{G_{t-s}(x-y)^2}{|z|^{2-2H}}\left[E \left(\int_0^s \int_{\bR^2}
|(G_{s-r}(y+z-v)-G_{s-r}(y-v))\sigma(u^n(r,v)) \right.\right.\\
& \; -\left.\left.(G_{s-r}(y+z-\bar v)-G_{s-r}(y-\bar v))\sigma(u^n(r,\bar v))|^2\frac{1}{|v-\bar v|^{2-2H}}\,dv \,d\bar v \, dr \right)^{p/2}\right]^{2/p} dy \, dz \,ds
\end{align*}
is uniformly bounded in $(t,x)$. After adding and subtracting the mixed term
$(G_{s-r}(y+z-v)-G_{s-r}(y-v))\sigma(u^n(r,\bar v))$ inside the squared term of the last integral above, we see that the previous integral is bounded by $2(A_1+A_2)$, where
\begin{align*}
A_1 &= \int_0^{t}ds \int_{\bR}dy \int_{\bR}dz\frac{G_{t-s}(x-y)^2}{|z|^{2-2H}}
\left[E\left(\int_0^s \int_{\bR^2} |G_{s-r}(y+z-v)-G_{s-r}(y-v)|^2 \right.\right.\\
&  \quad \qquad \times \left.\left.|\sigma(u^n(r,\bar v))-\sigma(u^n(r,\bar v))|^2 \frac{1}{|v-\bar v|^{2-2H}}\,dv \,d\bar v \, dr\right)^{p/2}\right]^{2/p} \\
A_2&= \int_0^{t}ds\int_{\bR} dy\int_{\bR} dz\frac{G_{t-s}(x-y)^2}{|z|^{2-2H}}\left[E \left(\int_0^s dr  \int_{\bR^2}\,dv \,d\bar v \, |\sigma(u^n(r,\bar v))|^2 \right. \right. \\
& \; \;\times \left. \left.\dfrac{|(G_{s-r}(y+z-v)-G_{s-r}(y-v))-(G_{s-r}(y+z-\bar v)-G_{s-r}(y-\bar v))|^2} {|v-\bar v|^{2-2H}} \right)^{p/2}\right]^{2/p}.
\end{align*}
We first deal with $A_2$. By Minkowski inequality,
\begin{align*}
A_2&\leq  \int_0^{t}ds\int_{\bR} dy\int_{\bR} dz\frac{G_{t-s}(x-y)^2}{|z|^{2-2H}}\left(\int_0^s dr  \int_{\bR^2}\,dv \,d\bar v \, \Big(E|\sigma(u^n(r,\bar v))|^p \Big)^{2/p}\right.  \\
&  \qquad \times \left.\dfrac{|(G_{s-r}(y+z-v)-G_{s-r}(y-v))-(G_{s-r}(y+z-\bar v)-G_{s-r}(y-\bar v))|^2} {|v-\bar v|^{2-2H}} \right).
\end{align*}
Taking into account that $\sigma$ is Lipschitz and $u^n$ satisfies
the second condition in \eqref{propertyP} (by the induction
hypothesis), we have that $A_2$ can be bounded by
\begin{align*}
& C \left[1+\left(\sup_{(t,x)}E|u^n(t,x)|^p\right)^{2/p}\right]\int_0^{t}ds\int_{\bR} dy\int_{\bR} dz\frac{G_{t-s}^2(x-y)}{|z|^{2-2H}}\\
&\; \times \int_0^s \int_{\bR^2}\dfrac{|(G_{s-r}(y+z-v)-G_{s-r}(y-v))-(G_{s-r}(y+z-\bar v)-G_{s-r}(y-\bar v))|^2} {|v-\bar v|^{2-2H}}\,dv \,d\bar v \, dr
\end{align*}
By Proposition \ref{propA}, all this expression can be estimated by (up to a constant)
$$\int_0^{t}ds\int_{\bR} dy\int_{\bR} dz\frac{G_{t-s}^2(x-y)}{|z|^{2-2H}}
\int_0^s \int_{\bR} |1-e^{-i\xi z}|^2 \, |\cF G_{s-r}(\xi)|^2 \,
|\xi|^{1-2H}d\xi \,dr.$$
Using Fubini's theorem and Lemma
\ref{elem-lem2}, this later expression is equal (up to a constant)
to
\begin{align*}
&\int_0^{t}\int_{\bR}G_{t-s}^2(x-y)
\left(\int_0^s \int_{\bR} |\cF G_{s-r}(\xi)|^2 \, |\xi|^{2(1-2H)}d\xi \,dr \right)\,dy \,ds, \\
& \leq \left(\int_0^{t}\int_{\bR}G_{t-s}^2(x-y)dy ds\right)
\left(\int_0^T \int_{\bR} |\cF G_r(\xi)|^2 \, |\xi|^{2(1-2H)}d\xi \,dr \right),
\end{align*}
which is uniformly bounded in $(t,x)$ by Lemma \ref{lemma-G1} with $\alpha=2(1-2H)$.

As far as $A_1$ is concerned, applying again Minkowski inequality
for integrals, we have:
\begin{align*}
A_1 &\leq \int_0^{t}ds \int_{\bR}dy \int_{\bR}dz\frac{G_{t-s}(x-y)^2}{|z|^{2-2H}}\left(\int_0^s \int_{\bR^2} |G_{s-r}(y+z-v)-G_{s-r}(y-v)|^2 \right.\\
&  \quad \qquad \times \left.\Big(E|\sigma(u^n(r,\bar v))-\sigma(u^n(r,\bar v))|^p \Big)^{2/p} \frac{1}{|v-\bar v|^{2-2H}}\,dv \,d\bar v \, dr\right)
\end{align*}
Using the fact that $\sigma$ is Lipschitz and Fubini's theorem,
\begin{multline*}
A_1 \leq C \int_0^{t}ds\int_0^{s} dr\int_{\bR} dv\int_{\bR} d\bar v
\dfrac{\Big(E|u^{n}(r,v)-u^n(r,\bar v)|^p\Big)^{2/p}}{|v-\bar v|^{2-2H}}\\
\times\Big( \int_{\bR} dz\int_{\bR} dy \dfrac{G_{t-s}(x-y)^2}{|z|^{2-2H}}
\,|G_{s-r}(y+z-v)-G_{s-r}(y-v)|^2\Big)
\end{multline*}
Observe that, doing the change of variables $\bar y=y-v$,
\begin{align*}
 & \int_{\bR} dz\intR dy \dfrac{G_{t-s}(x-y)^2}{|z|^{2-2H}}
\,|G_{s-r}(y+z-v)-G_{s-r}(y-v)|^2\\
& \quad = \intR dz\intR d\bar y \,\dfrac{G_{t-s}(x-\bar
y-v)^2}{|z|^{2-2H}} \,|G_{s-r}(\bar y+z)-G_{s-r}(\bar y)|^2.
\end{align*}
Hence, again by Fubini theorem,
\begin{align*}
A_1 & \leq C \inte{t}ds\inte{s}dr\intR dz\intR d\bar y
\frac{|G_{s-r}(\bar y+z)-G_{s-r}(\bar y)|^2}{|z|^{2-2H}}\\
& \quad \quad \times\Big(\intR dv\intR d\bar v \,G_{t-s}(x-\bar y-v)^2\,
\dfrac{\Big(E|u^{n}(r,v)-u^n(r,\bar v)|^p\Big)^{2/p}}{|v-\bar v|^{2-2H}}\Big)\\
& =\inte{t}ds\inte{s}d\bar r\intR dz\intR d\bar y
\frac{|G_{\bar r}(\bar y+z)-G_{\bar r}(\bar y)|^2}{|z|^{2-2H}}\\
& \quad \quad \times\Big(\intR dv\intR d\bar v \,G_{t-s}(x-\bar y-v)^2\,
\dfrac{\Big(E|u^{n}(s-\bar r,v)-u^n(s-\bar r,\bar v)|^p\Big)^{2/p}}{|v-\bar
v|^{2-2H}}\Big),
\end{align*}
where in the last equality we have done the change of variable $\bar
r=s-r$. By applying Fubini's theorem one more time, this last term is equal to
\begin{multline*}
\inte{t} d\bar r\intR dz\intR d\bar y \frac{|G_{\bar r}(\bar y+
z)-G_{\bar r}(\bar y)|^2}{|z|^{2-2H}}\\
 \times \int_{\bar r}^t ds\intR dv\intR d\bar v\, G_{t-s}(x-\bar
 y-v)^2\,\dfrac{\Big(E|u^{n}(s-\bar r,v)-u^n(s-\bar r,\bar v)|^p\Big)^{2/p}}{|v-\bar
v|^{2-2H}}.
\end{multline*}
Performing now the change of variables $\bar s=s-\bar r$, we can
write
\begin{multline*}
A_1 \leq C \inte{t} d\bar r\intR dz\intR d\bar y \frac{|G_{\bar
r}(\bar y+z)-G_{\bar r}(\bar y)|^2}{|z|^{2-2H}}\\
 \times \Big( \int_0^{t-\bar r} d\bar s\intR dv\intR d\bar v\, G_{t-\bar r-\bar s}(x-\bar y-v)^2\,\dfrac{\Big(E|u^{n}(\bar s,v)-u^n(\bar s,\bar v)|^p\Big)^{2/p}}{|v-\bar
v|^{2-2H}}\Big)\\
\le C \sup_{(\nu,w)\in[0,T]\times \bR} \left\{ \int_0^{\nu} d\bar
s\intR dv\intR d\bar v\, G_{\nu-\bar s}(
 w-v)^2\,\dfrac{\Big(E|u^{n}(\bar s,v)-u^n(\bar s,\bar v)|^p\Big)^{2/p}}{|v-\bar
v|^{2-2H}}\right\}\\
\times \inte{t} d\bar r\intR dz\intR d\bar y \frac{|G_{\bar
r}(\bar y+z)-G_{\bar r}(\bar y)|^2}{|z|^{2-2H}}.
\end{multline*}
By the induction hypothesis (the third condition in \eqref{propertyP}), the supremum appearing in the last term above is
finite. The remaining integral is
\begin{align*}
& \inte{t} d\bar r\intR dz\intR d\bar y  \frac{|G_{\bar
r}(\bar y+z)-G_{\bar r}(\bar y)|^2}{|z|^{2-2H}} =\frac{c_H}{C_H}\inte{t} d\bar r\intR |\cF G_{\bar r}(\xi)|^2\,|\xi|^{1-2H}\,d\xi,
\end{align*}
(by Proposition \ref{propA}), and this is uniformly bounded in $t\in
[0,T]$,  by Lemma \ref{lemma-G1}. This concludes the proof of
\eqref{bound2-n+1}.

\smallskip

\noindent {\it Step 4.} This final step is devoted to prove that
property \eqref{propertyQ} hold for $u^{n+1}$.

We consider first the space increments of $u^{n+1}$. By definition
\eqref{def-Picard} of $u^{n+1}$ and Theorem \ref{theoremA}, we see
that
$$E|u^{n+1}(t,x+h)-u^{n+1}(t,x)|^2 \leq C(I_0+I_1+I_2),$$
where $I_0=|w(t,x+h)-w(t,x)|^2$,
\begin{align*}
I_1 &= E\int_0^t \int_{\bR^2} |G_{t-s}(x+h-y)-G_{t-s}(x-y)|^2 \dfrac{|\sigma(u^n(s,y))-\sigma(u^n(s,z))|^2}{|y-z|^{2-2H}}\,dy \,dz\,ds\\
I_2 &= E \int_0^t \int_{\bR^2} \frac{|\sigma(u^n(s,z))|^2}{|y-z|^{2-2H}}\,
|(G_{t-s}(x+h-y)-G_{t-s}(x-y))-\\
& \qquad \qquad \qquad \qquad \qquad (G_{t-s}(x+h-z)-G_{t-s}(x-z))|^2\,dy \,dz \,ds.
\end{align*}
We have already proved that $I_0 \leq C|h|^{2H}$. Let us treat
$I_1$. Since $\sigma$ is Lipschitz,
\begin{align*}
& I_1 \leq  C E\int_0^t \int_{\bR} |G_{t-s}(x+h-y)-G_{t-s}(x-y)|^2 \left(\int_{\bR}\dfrac{|u^n(s,y+z)-u^n(s,y)|^2}{|z|^{2-2H}}\,dz\right) \,dy\,ds \\
& \quad =  C(I_1'+I_1''),
\end{align*}
where $I_1'$ and $I_1''$ denote  the integrals corresponding to the
regions $\{|z|>1\}$, respectively $\{|z| \leq 1\}$. Since
$\int_{\{|z|>1\}}|z|^{2H-2}dz$ is finite, we have
\begin{eqnarray*}
I_1' & \leq & C \sup_{(t,x) \in [0,T] \in \bR}E|u_n(t,x)|^2 \, \int_0^t \int_{\bR} |G_{t-s}(x+h-y)-G_{t-s}(x-y)|^2  \,dy\,ds\\
&=&  C \sup_{(t,x) \in [0,T] \in \bR}E|u^n(t,x)|^2 \int_0^t \int_{\bR} |1-e^{-i \xi h}|^2 \, |\cF G_s(\xi)|^2 \, d\xi \,ds \\
&=& C \sup_{(t,x) \in [0,T] \in \bR}E|u^n(t,x)|^2 \int_0^t \int_{\bR} (1-\cos(\xi h)) \, |\cF G_s(\xi)|^2 \, d\xi \,ds .
\end{eqnarray*}
using Plancherel's theorem for the first equality above.  Using
Lemma \ref{lemma-G3} (with $\alpha=0$), we obtain that
$$I_1' \leq C |h| \, \sup_{(t,x) \in [0,T] \in \bR} E|u^n(t,x)|^2.$$

For $I_1''$, we use the induction hypothesis (the first condition in \eqref{propertyQ}), to infer that $|u^n(s,y+z)-u^n(s,y)|^2 \leq C_n |z|^{2H}$.
Since $\int_{\{|z| \leq 1\}}|z|^{4H-2}dz$ is finite,
\begin{eqnarray*}
I_1'' & \leq & C C_n \int_0^t \int_{\bR} |G_{t-s}(x+y+h)-G_{t-s}(x+y)|^2 \, dy\, ds \leq C C_n |h|.
\end{eqnarray*}

Let us now treat the term $I_2$. Using the fact that $\sigma$ is
Lipschitz and Proposition \ref{propA},
\begin{eqnarray*}
I_2 & \leq & C \left(1+\sup_{(t,x)}E|u^n(t,x)|^2\right)\int_0^t \int_{\bR}|e^{-i \xi h}-1|^2 |\cF G_{t-s}(\xi)|^2 \,|\xi|^{1-2H}\,d\xi \,ds\\
&=&C \left(1+\sup_{(t,x)}E|u^n(t,x)|^2\right) \int_0^t \int_{\bR}(1-\cos(h|\xi|)) \, |\cF G_{t-s}(\xi)|^2 \,|\xi|^{1-2H}\,d\xi \,ds.
\end{eqnarray*}
Using Lemma \ref{lemma-G3} (with $\alpha=1-2H$), we obtain that
$$I_2 \leq C \left(1+\sup_{(t,x)}E|u^n(t,x)|^2\right)|h|^{2H}.$$
This concludes the proof for the space increments of  $u^{n+1}$.

We consider now the time increments of $u^{n+1}$. We assume that $h>0$. The case $h<0$ is similar. By \eqref{def-Picard} and Theorem \ref{theoremA}, for any $t<T-h$,
$$E|u^{n+1}(t+h,x)-u^{n+1}(t,x)|^2 \leq C(J_0+J_1+J_2),$$
where $J_0=|w(t+h,x)-w(t,x)|^2$ ,
\begin{eqnarray*}
J_1 &=& E\int_{t}^{t+h}\int_{\bR^2}\dfrac{|G_{t+h-s}(x-y)\sigma(u^n(s,y))-G_{t+h-s}(x-z)
\sigma(u^n(s,z))|^2}{|y-z|^{2-2H}}\,dy\,dz\,ds\\
J_2&=& E \int_0^t \int_{\bR^2}|(G_{t+h-s}(x-y)-G_{t-s}(x-y))\sigma(u^n(s,y))-\\
& & \quad \quad \qquad
(G_{t+h-s}(x-z)-G_{t-s}(x-z))\sigma(u^n(s,z))|^2|y-z|^{2H-2}\,dy\,dz \,ds.
\end{eqnarray*}
It was shown above that $J_0 \leq C|h|^{2H}$ (the case $n=0$). As
far as  $J_1$ is concerned, adding and subtracting
$G_{t+h-s}(x-y)\sigma(u^n(s,z))$, we see that $J_1 \leq
2(J_{11}+J_{12})$, where
\begin{eqnarray*}
J_{11} &=& E\int_t^{t+h} \int_{\bR^2} G_{t+h-s}^2(x-y) \dfrac{|\sigma(u^n(s,y))-\sigma(u^n(s,z))|^2}{|y-z|^{2-2H}}\,dz\,dy\,ds\\
J_{12} &=& E\int_{t}^{t+h} \int_{\bR^2}|\sigma(u^n(s,z))|^2 \,\dfrac{|G_{t+h-s}(x-y) -G_{t+h-s}(x-z)|^2}{ |y-z|^{2-2H}}\,dy\,dz \,ds.
\end{eqnarray*}
Since $\sigma$ is Lipschitz, $J_{11}$ is smaller than a constant
times
$$E\int_t^{t+h} \int_{\bR} G_{t+h-s}^2(x-y) \left(\int_{\bR}\dfrac{|u_n(s,y+z)-u_n(s,z)|^2}{|z|^{2-2H}}\,dz\right)\,dy\,ds,$$
which can be written $J_{11}'+J_{11}''$, where $J_{11}'$ and $J_{11}''$ are the integrals corresponding to the regions $\{|z|>1\}$, respectively $\{|z|\leq 1\}$. Since $H<1/2$, the integral $\int_{\{|z|>1\}}|z|^{2H-2}dz$ is finite, and hence
\begin{eqnarray*}
J_{11}' & \leq & C \sup_{(t,x)}E|u^n(t,x)|^2 \int_{t}^{t+h}\int_{\bR}G_{t+h-s}^2(x-y) \,dy\,ds\\
&=& C \sup_{(t,x)}E|u^n(t,x)|^2  \int_0^h \int_{\bR}G_s^2(s,y)\,dy \,ds.
\end{eqnarray*}
The last integral is equal to $Ch^2$ for the wave equation, respectively $Ch^{1/2}$ for the heat equation (see \eqref{L2normG}). For $J_{11}''$, we use the induction hypothesis (the first condition in \eqref{propertyQ}). Since $H>1/4$, $\int_{\{|z| \leq 1\}}|z|^{4H-2}dz$ is finite and
$$J_{11}'' \leq C C_n \int_{t}^{t+h}\int_{\bR}G_{t+h-s}^2(x-y) \,dy\,ds,$$
which is the same integral as above. As for $J_{12}$, since $\sigma$ is Lipschitz,
\begin{eqnarray*}
J_{12} & \leq & \left(1+\sup_{(t,x)}E|u^n(t,x)|^2 \right) \int_{t}^{t+h} \int_{\bR^2} \dfrac{|G_{t+h-s}(x-y)-G_{t+h-s}(x-z)|^2}{|y-z|^{2-2H}}\,dy\,dz \\
&=& C \left(1+\sup_{(t,x)}E|u^n(t,x)|^2 \right) \int_{0}^{h} \int_{\bR} |\cF G_{r}(\xi)|^2 \,|\xi|^{1-2H}\, d\xi dr,
\end{eqnarray*}
using Proposition \ref{propA} (after the change of variables $r=t+h-s$). By applying Lemma \ref{lemma-G1} with $\alpha=1-2H$, we see that the last integral is equal to $C h^{2H+1}$ for the wave equation, and $C h^{H}$ for the heat equation.

 We now treat $J_2$. As in the case of $J_1$, we add and subtract
the mixed term $(G_{t+h-s}(x-y)-G_{t-s}(x-y))\sigma(u^n(s,z))$.
Hence $J_2 \leq 2(J_{21}+J_{22})$, where
\begin{eqnarray*}
J_{21} &=& E \int_0^t \int_{\bR^2}|G_{t+h-s}(x-y)-G_{t-s}(x-y)|^2\\
& & \quad \quad \qquad |\sigma(u^n(s,y))-\sigma(u^n(s,z))|^2 \,|y-z|^{2H-2}\,dy \,dz \,ds\\
J_{22} &=& E \int_0^t \int_{\bR^2} |\sigma(u^n(s,z))|^2 \, |(G_{t+h-s}(x-y)-G_{t-s}(x-y))-\\
& & \quad \quad \qquad (G_{t+h-s}(x-z)-G_{t-s}(x-z))|^2 \,|y-z|^{2H-2}\,dy\,dz\,ds.
\end{eqnarray*}
These terms are treated similarly to $J_{11}$, respectively $J_{12}$. More precisely, $J_{21} \leq C(J_{21}'+J_{21}'')$, where $J_{21}'$ and $J_{21}''$ are integrals corresponding to the regions $\{|z|>1\}$, respectively $\{|z| \leq 1\}$. Similarly to $J_{11}'$, we obtain that
$$J_{21}' \leq  C \sup_{(t,x)}E|u^n(t,x)|^2 \int_0^t \int_{\bR}|G_{t+h-s}(x-y)-G_{t-s}(x-y)|^2\,dy\,ds.$$
By Plancherel's theorem, the previous integral is equal to
$(2\pi)^{-1}$ times
$$I:=\int_{0}^{t} \int_{\bR}|\cF G_{s+h}(\xi)-\cF G_{s}(\xi)|^2 \, d\xi \,ds.$$
By applying Lemma \ref{lemma-G4} with $\alpha=0$, we see that this integral is smaller than $C t h$ for the wave equation, respectively $C h^{1/2}$ for the heat equation.

Similarly to $J_{11}''$, we have
$$J_{21}'' \leq C C_n \int_0^t \int_{\bR}|G_{t+h-s}(x-y)-G_{t-s}(x-y)|^2 \,dy \,ds.$$
As noted above, this last integral is bounded by $Cth$ for the wave equation, respectively $C h^{1/2}$ for the heat equation.

As for $J_{22}$, similarly to the argument used for $J_{12}$, we have
\begin{eqnarray*}
& & J_{22}  \leq  C\left(1+\sup_{(t,x)}E|u^n(t,x)|^2\right)\int_0^t \int_{\bR^2}|(G_{t+h-s}(x-y)-G_{t-s}(x-y))-\\
& & \qquad \quad \qquad \qquad (G_{t+h-s}(x-z)- G_{t-s}(x-z))|^2 \,|y-z|^{2H-2}\,dy\,dz\,ds.
\end{eqnarray*}
By Proposition \ref{propA}, the integral above is equal to a constant times
$$\int_0^t \int_{\bR}|\cF G_{t+h-s}(\xi)-\cF G_{t-s}(\xi)|^2 \, |\xi|^{1-2H}\,d\xi\,ds$$
By applying Lemma \ref{lemma-G4} with $\alpha=1-2H$, we see that
this integral is smaller than $C t h^{2H}$ for the wave equation,
respectively $Ch^{H}$ for the heat equation. This concludes the
proof for the space increments of $u^{n+1}$. To summarize, we have
proved that $u^{n+1}$ satisfies the property \eqref{propertyQ} with
$C_{n+1}=C(1+K_n+C_n)$, where $K_n=\sup_{(t,x) \in [0,T] \in \bR}
E|u^n(t,x)|^2$.

Therefore, we can conclude the proof of Theorem
\ref{Picard-well-def}. \qed

\subsection{Convergence of Picard iterations}
\label{section-conv-Picard}

In this section, we will prove that,
for any $(t,x)\in [0,T]\times \bR$, the sequence $\{u^n(t,x),\, n\geq 0\}$ of random variables
converges in $L^p(\Omega)$, for $p \geq 2$ fixed. We will assume that $\sigma$ satisfies the condition:
\begin{equation}
\label{sigma-affine}
|\sigma(x)-\sigma(y)-\sigma(u)+\sigma(v)| \leq C|x-y-u+v|
\end{equation}
for any $x,y,u,v \in \bR$, for some constant $C>0$. Note that \eqref{sigma-affine} holds if and only if $\sigma$ is affine.
For any $n \geq 1$ and for any $t\in [0,T]$, we define
\begin{align*}
V_{n}(t)& :=\sup_{x\in \bR}\Big(E |u^n(t,x)-u^{n-1}(t,x)|^p \Big)^{2/p}\\
 W_n(t) & := \sup_{x\in \bR}  \int_0^t \!\! \int_{\bR^2} G_{t-s}^2(x-y)\,|y-z|^{2H-2}\\
 & \times
 \Big(E|u^n(s,y)-u^{n-1}(s,y) - u^n(s,z)+u^{n-1}(s,z)|^p\Big)^{2/p} \, dy\,dz\,ds.
\end{align*}


The following result establishes a recurrence relation for the pair $(V_n,W_n)$.

\begin{theorem}
\label{theorem-VW}
Assume that $\sigma$ is an affine function, i.e. $\sigma(x)=ax+b$ for some $a,b \in \bR$. For any $n \geq 0$ and for any $t\in [0,T]$, we have
\begin{equation}
\label{VW-rec1}
V_{n+1}(t)\leq \int_0^t V_{n}(s) J_1(t-s) ds + C\, W_{n}(t),
\end{equation}
and
\begin{equation}
\label{VW-rec2}
W_{n+1}(t)\leq \int_0^t V_{n}(s) J_2(t-s) ds + \int_0^t W_{n}(s) J_1(t-s) ds,
\end{equation}
where $C=z_p^{p/2}C_H 2^{2(p-1)/p} a^2$, and $J_1$ and $J_2$ are some non-negative functions in $L^1([0,T])$.
\end{theorem}

\noindent {\bf Proof:} As usual, we denote by $C$ a constant (depending on $p,H,T$ and $a$) which may be different from line to line.
We split the proof in two parts: the first one will be devoted to prove \eqref{VW-rec1}, while in the second
one we will show \eqref{VW-rec2}.

\smallskip

\noindent {\it Step 1.}
We start by checking (\ref{VW-rec1}). By \eqref{def-Picard} and Theorem \ref{theoremA},
\begin{align*}
& E |u^{n+1}(t,x)-u^{n}(t,x)|^p \\
 & \quad =
 E  \left|\int_0^t\!\!\int_{\bR} G_{t-s}(x-y)\big\{\sigma(u^{n}(s,y))-\sigma(u^{n-1}(s,y))\big\} X(ds,dy)
 \right|^p \\
 & \quad \leq C E \left(\int_0^t\!\!\int_{\bR^2} \big| G_{t-s}(x-y)\big\{\sigma(u^{n}(s,y))-\sigma(u^{n-1}(s,y))\big\} \right.\\
 & \left.\quad \qquad - G_{t-s}(x-z)\big\{\sigma(u^{n}(s,z))-\sigma(u^{n-1}(s,z))\big\}\big|^2 |y-z|^{2H-2} \, dy\, dz\,ds \right)^{p/2}.
\end{align*}
After adding and subtracting the term
$$G_{t-s}(x-y)\big\{\sigma(u^{n}(s,z))-\sigma(u^{n-1}(s,z))\big\},$$
we see that $\Big(E |u^{n+1}(t,x)-u^{n}(t,x)|^p\Big)^{2/p} \leq C(A_1+A_2)$, where
\begin{align*}
A_1 & =\left[E \left(\int_0^t\!\!\int_{\bR^2}  G_{t-s}^2(x-y) \,  |y-z|^{2H-2}\right.\right.\\
 &  \;\left.\left.\times
 \big| \sigma(u^{n}(s,y))-\sigma(u^{n-1}(s,y)) - \sigma(u^{n}(s,z)) + \sigma(u^{n-1}(s,z))   \big|^2 \,dy\,dz\,ds\right)^{p/2} \right]^{2/p}\\
 A_2  & = \left[E \left(\int_0^t\!\!\int_{\bR^2}  |G_{t-s}(x-y) - G_{t-s}(x-z)|^2 \right. \right.\\
 & \left. \left. \qquad \quad \times  \big|\sigma(u^{n}(s,z))-\sigma(u^{n-1}(s,z))\big|^2 |y-z|^{2H-2} dydzds \right)^{p/2} \right]^{2/p}.
\end{align*}
As far as $A_1$ is concerned, applying \eqref{sigma-affine} (with $C=a$) and Minkowski's inequality for integrals,
we obtain that
\begin{equation}
 A_1 \leq a^2 \, W_{n}(t).
 \label{eq:998}
\end{equation}
Regarding $A_2$, using the fact that $\sigma(x)=ax+b$ and applying
again Minkowski's inequality for integrals, we have
\begin{align}
 A_2 & \leq a^2  \int_0^t\!\!\int_{\bR^2}  |G_{t-s}(x-y) - G_{t-s}(x-z)|^2    \nonumber\\
 & \quad \qquad \times \Big(E|u^n(s,z)-u^{n-1}(s,z)|^p \Big)^{2/p}\, |y-z|^{2H-2}\,dy\,dz\,ds \nonumber \\
 & \leq a^2   \int_0^t V_{n}(s) \int_{\bR^2} |G_{t-s}(x-y) - G_{t-s}(x-z)|^2 |y-z|^{2H-2} dydz ds  \nonumber \\
 & =\int_0^t V_{n}(s) J_1(t-s)ds,
 \label{eq:997}
\end{align}
where
\begin{equation}
\label{def-J1}
 J_1(t) := a^2 \int_{\bR^2} |G_{t}(y) - G_{t}(z)|^2 |y-z|^{2H-2} dydz  = C \int_{\bR} |\cF G_{t}(\xi)|^2 |\xi|^{1-2H} d\xi,
\end{equation}
by Proposition \ref{propA}. By \eqref{int-wave} and \eqref{int-heat}, $J_1(t)=C t^{2H}$ for the wave equation and $J_1(t)=C t^{H-1}$ for the heat equation.
In both cases, $J_1$ is an integrable function on $[0,T]$.
Since estimates (\ref{eq:998}) and (\ref{eq:997}) are uniform with respect to $x\in \bR$, putting them together
one obtains that (\ref{VW-rec1}) holds.

\medskip

\noindent {\it Step 2.}
 Let us now prove (\ref{VW-rec2}). For this, we first define, for any $n\geq 1$,
$$m_n(r,v):= \sigma(u^n(r,v))-\sigma(u^{n-1}(r,v)), \quad (r,v)\in [0,T]\times \bR.$$
Taking into account that $W_{n+1}(t)$ can be written as
\begin{align*}
 W_{n+1}(t) & = \sup_{x\in \bR} \, \int_0^t \!\! \int_{\bR^2} \frac{G_{t-s}(x-y)^2}{|z|^{2-2H}} \\
 & \quad \times
 \Big(E |u^{n+1}(s,y+z)-u^{n}(s,y+z) - u^{n+1}(s,y)+u^{n}(s,y)|^p \Big)^{2/p} \, dy\,dz\,ds,
\end{align*}
we see that the latter integral is bounded (up to a constant) by
\begin{align}
\nonumber
& \int_0^t \int_{\bR^d}\dfrac{G_{t-s}^2(x-y)}{|z|^{2-2H}} \left[E \left(\int_0^s \int_{\bR^2}|(G_{s-r}(y+z-v)-G_{s-r}(y-v))m_n(r,v) \right. \right.\\
\label{W-int}
& \quad - \left. \left. (G_{s-r}(y+z-\bar v)-G_{s-r}(y-\bar v))m_n(r,\bar v)|^2\frac{1}{|v-\bar v|^{2-2H}} dv \, d\bar v \,dr \right)^{p/2} \right]^{2/p}dy\,dz\,ds
\end{align}
This follows by Theorem \ref{theoremA}, noting that (by \eqref{def-Picard})
\begin{align*}
& E|u^{n+1}(s,y+z)-u^n(s,y+z) - u^{n+1}(s,y)+u^n(s,y)|^p\\
& \qquad =E\left|\int_0^s \int_{\bR}(G_{s-r}(y+z-v)-G_{s-r}(y-v))m_n(r,v)X(dr,dv)\right|^p.
\end{align*}
Adding and subtracting the term $\{G_{s-r}(y+z-v)-G_{s-r}(y-v)\} m_{n}(r,\bar v)$, we see that the integral \eqref{W-int} is bounded by $2(B_1+B_2)$, where
\begin{align*}
 B_1 & = \int_0^t ds \int_{\bR}dy \int_{\bR}dz \frac{G_{t-s}(x-y)^2}{|z|^{2-2H}}
 \left[E\int_0^s \left(\!\! \int_{\bR^2} |G_{s-r}(y+z-v)-G_{s-r}(y-v)|^2 \right. \right.\\
 & \quad  \quad \left. \left. \times |m_{n}(r,v)-m_{n}(r,\bar v)|^2 \frac{1}{|v-\bar v|^{2-2H}} dv d\bar v dr \right)^{p/2} \right]^{2/p} \\
 B_2
 &  = \int_0^t ds\int_{\bR}dy\int_{\bR}dz \frac{G_{t-s}(x-y)^2}{|z|^{2-2H}}
 \left[E\left(\int_0^sdr\int_{\bR^2}dv d\bar v\, |m_{n}(r,\bar v)|^2 \right.\right.\\
 &  \quad  \times \left. \left. \dfrac{|(G_{s-r}(y+z-v)-G_{s-r}(y-v)) - (G_{s-r}(y+z-\bar v) - G_{s-r}(y-\bar v))|^2}{|v-\bar v|^{2-2H}}  \right)^{p/2} \right]^{2/p}.
\end{align*}

To deal with the term $B_1$, we argue as we did for the term $A_1$ in
the proof of Theorem \ref{Picard-well-def}. Indeed, by applying Minkowski's inequality for integrals,
 \begin{align*}
B_1 & \leq \int_0^t ds \int_{\bR}dy \int_{\bR}dz \frac{G_{t-s}(x-y)^2}{|z|^{2-2H}}
 \int_0^s  \int_{\bR^2} |G_{s-r}(y+z-v)-G_{s-r}(y-v)|^2 \\
 & \quad  \quad  \times \Big(E|m_{n}(r,v)-m_{n}(r,\bar v)|^p \Big)^{2/p} \frac{1}{|v-\bar v|^{2-2H}} dv d\bar v dr
 \end{align*}
We now apply several times Fubini's theorem, together with the following changes of variables:
\begin{enumerate}
 \item $\bar y:=y-v$ ($v$ fixed),
 \item $\bar r:= s-r$ ($s$ fixed),
 \item $\bar s:=s-\bar r$ ($\bar r$ fixed).
\end{enumerate}
Using these techniques, we infer that
\begin{align*}
B_1 &\leq C \int_0^t d\bar{r}\int_{\bR}dz \int_{\bR}d\bar y \, \dfrac{|G_{\bar r}(\bar y+z)-G_{\bar r}(\bar y)|^2}{|z|^{2-2H}} \\
 & \qquad \times \left( \int_0^{t-\bar r}d\bar s \int_{\bR}dv \int_{\bR}d \bar v
\, G_{t-\bar r-\bar s}(x-\bar y-v)^2
 \dfrac{\Big(E|m_{n}(\bar s,v)- m_{n}(\bar s,\bar v)|^p\Big)^{2/p}}{|v-\bar v|^{2-2H}} \right)
\end{align*}
Using the fact that $\sigma$ satisfies \eqref{sigma-affine} and taking the supremum  with respect to $x-\bar y$ of the integral in the parenthesis above, we get
$$B_1 \leq C \int_0^t W_{n}(t-\bar r) \int_{\bR^2} \frac{|G_{\bar r}(y+\bar z)-G_{\bar r}(\bar y)|^2}{|z|^{2-2H}}
  d\bar y d z d\bar r.$$
Performing again another change of coordinates ($\tau:=t-\bar r$), and recalling definition \eqref{def-J1} of $J_1$, we finally obtain that
\begin{equation}
B_1\leq \int_0^t W_{n}(\tau) J_1(t-\tau) d\tau.
\label{eq:995}
\end{equation}

Concerning the term $B_2$, analogous arguments as before yield
\begin{align*}
B_2 & \leq \int_0^t ds\int_{\bR}dy\int_{\bR}dz \frac{G_{t-s}(x-y)^2}{|z|^{2-2H}}
 \int_0^s dr\int_{\bR^2}dv d\bar v\, \Big(E|m_{n}(r,\bar v)|^p \Big)^{2/p} \\
 & \quad \times  \dfrac{|(G_{s-r}(y+z-v)-G_{s-r}(y-v)) - (G_{s-r}(y+z-\bar v) - G_{s-r}(y-\bar v))|^2}{|v-\bar v|^{2-2H}}  \\
 & \leq C \int_0^t ds\int_{\bR}dy \int_{\bR}dz  \frac{G_{t-s}^2(x-y)}{|z|^{2-2H}} \int_0^s dr V_{n}(r)\\
 & \quad \times  \left(\int_{\bR^2}
 \frac{|G_{s-r}(y+z-v)-G_{s-r}(y-v) - G_{s-r}(y+z-\bar v) + G_{s-r}(y-\bar v)|^2}{|v-\bar v|^{2-2H}}
 dv d\bar v \right)
\end{align*}
By Proposition \ref{propA}, the expression inside the parenthesis above is  equal to
$$C
 \int_{\bR}|1-e^{-i\xi z}|^2\, |\cF G_{s-r}(\xi)|^2\,|\xi|^{1-2H}d\xi.$$
We apply Fubini's theorem and then Lemma \ref{elem-lem2} to compute the $dz$ integral. We end up with
$$B_2 \leq C \int_0^t\!\!\int_{\bR}  G_{t-s}^2(x-y)  \int_0^s V_{n}(r)
 \int_{\bR} |\cF G_{s-r}(\xi)|^2|\xi|^{2(1-2H)}d\xi dr dy ds.$$
Applying Fubini's theorem one more time, we obtain
\begin{align*}
 B_2 & \leq  \int_0^t V_{n}(r) \left[ \int_{r}^t\!\!\int_{\bR} G_{t-s}^2(x-y)
 \int_{\bR} |\cF G_{s-r} (\xi)|^2 |\xi|^{2(1-2H)}d\xi dy ds\right] dr \\
 & = \int_0^t V_{n}(r) \left[ \int_{r}^t\!\!\int_{\bR} G_{t-s}^2(z)
 \int_{\bR} |\cF G_{s-r} (\xi)|^2 |\xi|^{2(1-2H)}d\xi dz ds\right] dr
\end{align*}
All that we need to check is that the latter expression inside the brackets can be written in the form
$J_2(t-r)$, where $J_2:[0,T]\rightarrow \bR_+$ is an integrable function.
This is precisely the statement in Lemma \ref{lemma-G2}. Therefore,
\begin{equation}
B_2\leq \int_0^t V_{n}(r) J_2(t-r) dr.
\label{eq:993}
\end{equation}
Putting together (\ref{eq:995}) and (\ref{eq:993}), we conclude that
(\ref{VW-rec2}) holds. \qed

\vspace{3mm}

We now introduce the solution space. Let $p \geq 2$ be fixed.
We denote by $\cX$ the set of $L^2(\Omega)$-continuous and adapted processes $Y=\{Y(t,x);t \in [0,T],x\in \bR\}$ such that
$\|Y\|_{\cX_1}<\infty$ and $\|Y\|_{\cX_2}<\infty$, where
$$\|Y\|_{\cX_1}=\sup_{(t,x) \in [0,T] \times \bR}\Big(E|Y(t,x)|^p\Big)^{1/p}$$
and
\begin{align*}
\|Y\|_{\cX_2}&=\sup_{(t,x) \in [0,T] \times \bR} \left(\int_0^t \int_{\bR^2}G_{t-s}^{2}(x-y)\dfrac{\Big(E|Y(s,y)-Y(s,z)|^p\Big)^{2/p}}{|y-z|^{2-2H}}\,dy\,dz\,ds \right)^{1/2}
\end{align*}
For any $Y \in \cX$, we define $\|Y\|_{\cX}:=\|Y\|_{\cX_1}+\|Y\|_{\cX_2}$.
We identify two processes $Y_1$ and $Y_2$ for which $\|Y_1-Y_2\|_{\cX}=0$.
With these definitions, one checks that $(\cX,\|\cdot\|_{\cX})$ defines a Banach space.

\vspace{3mm}

We are now ready to prove the main result of this section.

\begin{theorem}
\label{main-th}
The sequence $(u^n)_{n \geq 0}$ converges in $\cX$ to a process $u$, which is $L^2(\Omega)$-continuous,
and is the unique solution in $\cX$ to equation \eqref{wave} (or \eqref{heat}).
\end{theorem}

\noindent {\bf Proof:} We first show the existence of $u$.
 By Theorem \ref{theorem-VW}, we have
$$M_{n+1}(t) \leq \int_0^t (M_n(s)+M_{n-1}(s))J(t-s)ds$$
where $J(t)=C(J_1(t)+J_2(t))$ and
\begin{align*}
M_n(t)&=\sup_{x\in \bR}\left\{\Big(E|u^n(t,x)-u^{n-1}(t,x)|^p\Big)^{2/p} \right\} +
\sup_{x\in \bR}\left\{ \int_0^t \int_{\bR^2}G_{t-s}^2(x-y)|y-z|^{2H-2}\right. \\
& \left. \qquad \times \Big(E|u^n(s,y)-u^{n-1}(s,y)-u^{n}(s,z)+u^{n-1}(s,z)|^p \Big)^{2/p} \,dy\,dz\,ds\right\}.
\end{align*}

By Lemma \ref{gronwall} below,
$$\sum_{n \geq 0}\sup_{t \in [0,T]}M_n(t)^{1/2}<\infty.$$
This implies that
$\sum_{n \geq 0}\|u^n-u^{n-1}\|_{\cX_i}<\infty$ for $i=1,2$, and consequently,
$\sum_{n \geq 0}\|u^n-u^{n-1}\|_{\cX}<\infty$. Therefore, $(u^n)_n$ is Cauchy in $\cX$.
Therefore, there exists a process $u \in \cX$ such that \begin{equation}
\label{un-to-u}
\lim_{n \to \infty}\|u^n-u\|_{\cX}=0.
\end{equation}
In particular, the process $u$ is $L^2(\Omega)$-continuous, adapted and satisfies
\[
 \sup_{(t,x) \in [0,T] \times \bR}\Big(E|u(t,x)|^p\Big)^{1/p} <\infty.
\]
By an extension to random fields of Proposition 3.21 of \cite{PZ07}, it follows that $u$ has a predictable modification (called also $u$).

We now show that $u$ satisfies the integral equation \eqref{int-eq}. For this, we take the limit as $n \to \infty$ in the definition \eqref{def-Picard} of the Picard sequence $(u^n)_n$. On the left-hand side, $u^{n+1}(t,x) \to u(t,x)$ in $L^p(\Omega)$ by \eqref{un-to-u}. The right-hand side of \eqref{def-Picard} converges in $L^p(\Omega)$ to the right-hand side of \eqref{int-eq}, since
\begin{equation}
\label{RHS-converges}
E\left|\int_0^t \int_{\bR}G_{t-s}(x-y)\{\sigma(u^n(s,y))-\sigma(u(s,y))\}X(ds,dy)\right|^p \to 0.
\end{equation}
To see this, we note that by Theorem \ref{theoremA}, the previous expectation is bounded by
\begin{align*}
& CE\left(\int_0^t \int_{\bR^2}|G_{t-s}(x-y)\{\sigma(u^n(s,y))-\sigma(u(s,y))\} \right.\\
& \left.\qquad \quad - G_{t-s}(x-z)\{\sigma(u^n(s,z))-\sigma(u(s,z))\}|^2 \,|y-z|^{2H-2}dydzds\right)^{p/2},
\end{align*}
which is bounded by $C(A_1+A_2)$, where
\begin{align*}
A_1 &= E\left(\int_0^t \int_{\bR^2}G_{t-s}^2(x-y) |y-z|^{2H-2} \right.\\
&\left.\qquad \times |\{\sigma(u^n(s,y))-\sigma(u(s,y))\}-\{\sigma(u^n(s,z))-\sigma(u(s,z))\}|^2 \,dy\,dz\,ds\right)^{p/2}\\
A_2 &=E\left(\int_0^t \int_{\bR^2}\dfrac{|G_{t-s}(x-y)-G_{t-s}(x-z)|^2}{|y-z|^{2-2H}} |\sigma(u^n(s,z))-\sigma(u(s,z))|^2 \,dy\,dz\,ds\right)^{p/2}
\end{align*}

Using Minkowski inequality and the fact that $\sigma$ satisfies \eqref{sigma-affine},
we see that $A_1 \leq C \|u^n-u\|_{\cX_2}^p$, which converges to $0$ by \eqref{un-to-u}.
As for $A_2$, using similar arguments we obtain
\begin{align*}
A_2 & \leq C \left(\sup_{(t,x)}E|u^n(t,x)-u(t,x)|^p\right)^{2/p} \\
& \qquad \; \times \left(\int_{0}^t \int_{\bR^2}\dfrac{|G_{t-s}(x-y)-G_{t-s}(x-z)|^2}{|y-z|^{2-2H}} \,dy\,dz\,ds\right)^{p/2},
\end{align*}

which also converges to $0$ by \eqref{un-to-u}. This concludes the proof of \eqref{RHS-converges}.

Finally, we prove the uniqueness of the solution in $\cX$. Assume that there exist two predictable processes $u$ and $v$ which both satisfy \eqref{int-eq}. We denote $d(t,x)=u(t,x)-v(t,x)$.
For any $t\in [0,T]$, we define
\begin{align*}
V(t)& :=\sup_{x\in \bR}E |d(t,x)|^2\\
W(t) & := \sup_{x\in \bR}  \int_0^t \!\! \int_{\bR^2} G_{t-s}^2(x-y)\,
 \frac{E|d(s,y) - d(s,z)|^2}{|y-z|^{2-2H}} \, dy\,dz\,ds.
\end{align*}
As in the proof of Theorem \ref{theorem-VW} (replacing $u^n$ by $u$ and $u^{n-1}$ by $v$, and taking $p=2$), we have:
\begin{align*}
V(t) &\leq \int_0^t V(s)J_1(t-s)+CW(t)\\
W(t) & \leq \int_0^t V(s)J_2(t-s)ds+\int_0^t W(s)J_1(t-s)ds.
\end{align*}
We denote $M(t)=V(t)+W(t)$. It follows that
$$M(t) \leq \int_0^t M(s)J(t-s)ds$$
where $J(t)=c(J_1(t)+J_2(t))$. By Lemma 15 of \cite{dalang99}, $M(t)=0$ for all
$t \in [0,T]$. Hence $u(t,x)=v(t,x)$ a.s. for any $t \in [0,T], x \in \bR$. \qed

\medskip

Note that our main result Theorem \ref{main-thm} is an immediate consequence of
the above Theorem \ref{main-th}.

\subsection{An extension of Gronwall's lemma}
\label{section-gronwall}

In this section, we state and prove an extension of Gronwall's lemma, which was used in the proof of Theorem \ref{main-th}. This result can be viewed as a version of Lemma 15 in \cite{dalang99}.

\begin{lemma}
\label{gronwall}
Let $(f_n)_{n \geq 0}$ be a sequence of non-negative functions on
$[0,T]$, such that $M_0=\sup_{t \in [0,T]}f_0(t)<\infty$ and $M_1=\sup_{t \in [0,T]}f_1(t)<\infty$. Let $M=M_0+M_1$. Assume that for any $n\ge 2$ and for any $t\in[0,T]$,
\begin{equation}\label{ineq-Gron}
f_{n}(t)\le\int_0^t \big(f_{n-1}(s)+f_{n-2}(s)\big)g(t-s)\,\, ds,
\end{equation}
where $g:[0,T]\rightarrow \mathbb R_+$ is an
integrable function. Then, there exists a sequence $(a_n)_{n\ge 1}$  of positive
numbers such that $\sum_{n \geq 0}a_n^{1/p}<\infty$ for any $p>0$, and
\begin{equation}
\label{fn-bounded-an}
\sup_{t\in[0,T]}f_n(t)\le M a_n \quad \mbox{for all} \ n \geq 0.
\end{equation}
In particular,
$\sum_{n \geq 0}\sup_{t \in [0,T]}f_n(t)^{1/p}<\infty$ for any $p>0$.
\end{lemma}

\noindent {\bf Proof:} The argument is similar to the one used in the proof of Lemma 15 of \cite{dalang99}. Set $G(t)=\int_0^t g(s) \,ds$. We assume that $G(T)>0$. (The case $G(T)=0$ is trivial.) Let $X,(X_n)_{n\ge 1}$ be
independent identically distributed random variables with values in $[0,T]$ and density $g(s)/G(T),s \in [0,T]$, defined on the same probability space $(\Omega,\cF,P)$. We denote $S_n=\sum_{i=1}^{n}X_i$.

Let $K=\max\big(G(T),1\big)$. Applying (\ref{ineq-Gron}) with $n=2$, we have:
\begin{align}\nonumber
f_2(t) & \le\int_0^t \Big(f_1(t-s)+f_0(t-s)\Big)g(s)ds\\
& =G(T)E \Big(1_{\{X_1\le t\}} \big(f_1(t-X_1)+f_0(t-X_1)\big)\Big)
\label{ineq-n-2}\\
\nonumber
& \le 2M\,K\, E(1_{\{X_1\le t\}}).
\end{align}
For $n=3$, we obtain
\begin{align*}
f_3(t) & \le\int_0^t \Big(f_2(t-s)+f_1(t-s)\Big)g(s)ds\\
& =G(T)\int_{\Omega} 1_{\{X(\omega)\le
t\}}\Big(f_2(t-X(\omega))+f_1(t-X(\omega))\Big)P(d\omega)\\
& \le G(T) \int_{\Omega} 1_{\{X(\omega)\le t\}}
\Big[G(T) \int_{\Omega}1_{\{X_1(\omega_1)\le
t-X(\omega)\}}\Big(f_1(t-X(\omega)-X_1(\omega_1))\\
& \qquad +f_0(t-X(\omega)-X_1(\omega_1))\Big)P(d\omega_1)
+f_1(t-X(\omega))\Big]P(d\omega)\\
& \leq G(T)\int_{\Omega}1_{\{X(\omega) \leq t\}} \Big[2M\,G(T) \int_{\Omega} 1_{\{X_1(\omega_1) \leq t-X(\omega)\}}P(d\omega_1)+M\Big]P(d\omega)\\
& \le2M \, G(T)^2 E(1_{\{X_1+X_2\le t\}})+M\,G(T)E(1_{\{X_1\le
t\}})\\
& \le 3M K^2 E(1_{\{X_1\le t\}}),
\end{align*}
where we used \eqref{ineq-n-2} for the second inequality above.
We denote $$k_n=\left\lfloor\frac{n}2\right\rfloor,$$ where
$\lfloor x\rfloor$ is the integer $k$ such that $k \leq x <k+1$.

Next we show that, for any $n\ge 2$, and for any $t \in [0,T]$, we have
\begin{equation}
\label{induc}
f_n(t)\le  M\,b_{n+1}K ^{n-1}P(S_{k_n}\le t),
\end{equation}
where $b_n$ is the $n$-th term of Fibonacci sequence with
$b_1=b_2=1$. The proof of (\ref{induc}) follows by induction. We
have already seen that it holds for $n=2$ and $n=3$. Suppose that
(\ref{induc}) holds up to some $n$. We prove it for $n+1$. We must
distinguish two cases: $n$ is odd and $n$ is even.

First, assume that $n$ is odd. Say $n=2m+1$ for some positive integer $m$. Then $k_n=k_{n-1}=m$ and $k_{n+1}=k_n+1=m+1$. Thus,
\begin{align*}
f_{n+1}(t) & \le\int_0^t \Big(f_n(t-s)+f_{n-1}(t-s)\Big)g(s)ds\\
& =G(T)\int_{\Omega}1_{\{X(\omega)\le t\}}\Big(f_n(t-X(\omega))+f_{n-1}(t-X(\omega))\Big)P(d\omega)\\
& \leq G(T) \int_{\Omega} 1_{\{X(\omega) \leq t\}} \Big(
M b_{n+1} K^{n-1} E(1_{\{S_{k_n} \leq t-X(\omega)\}}) \\
& \qquad  \quad + M b_n K^{n-2} E(1_{\{ S_{k_{n-1}} \leq t-X(\omega)\}})
\Big) P(d\omega),
\end{align*}
using the induction hypothesis for the last inequality.
Using the fact that $k_{n-1}=k_{n}$ and the recurrence relation of Fibonacci
numbers:
$$b_{n+2}=b_{n+1}+b_{n},$$
we obtain that the right-hand side of the last inequality
is bounded by
\begin{align*}
& G(T)M b_{n+2} K^{n-1} \int_{\Omega}1_{\{X(\omega) \leq t\}}E[1_{\{S_{k_n} \leq t-X(\omega)\}}]P(d\omega) \\
& = G(T) M b_{n+2} K^{n-1} \int_{\Omega} 1_{\{X_m(\omega_1) \leq t\}} \int_{\Omega}1_{\{S_m(\omega_2) \leq t-X_m(\omega_1)\}}P(d\omega_2) P(d\omega_1)\\
& \leq M b_{n+2} K^{n}\,P(S_{m+1} \leq t)=M b_{n+2} K^{n}\,P(S_{k_{n+1}} \leq t).
\end{align*}

Suppose now that $n$ is even. Say $n=2m$ for a positive integer $m$. Then
$k_n=k_{n+1}=m$ and $k_{n-1}=k_n-1=m-1$. Using again the induction hypothesis, we have
\begin{align*}
f_{n+1}(t) & \le G(T)\int_{\Omega} 1_{\{X(\omega)\le
t\}}\Big(f_n(t-X(\omega))+f_{n-1}(t-X(\omega))\Big)P(d\omega)\\
& \le G(T) \int_{\Omega} 1_{\{X(\omega)\le t\}} \Big(M\, b_{n+1}\,K^{n-1}
E(1_{\{S_{k_n}\le t-X(\omega)\}})+\\
& \quad \qquad M \,b_{n}\,K^{n-2}E(1_{\{S_{k_{n-1}} \le
t-X(\omega)\}}) \Big)P(d\omega)\\
& \leq M b_{n+1}K^{n}
P(S_{k_n} \leq t)+\\
& \qquad \quad M b_n K^{n-1}\int_{\Omega}1_{\{X_m(\omega_1) \leq t\}}\int_{\Omega} 1_{\{S_{m-1}(\omega_2) \leq t-X_m(\omega_1)\}}P(d\omega_2)P(d\omega_1)\\
& \leq  M b_{n+1}K^{n} P(S_{k_n} \leq t)+ M b_n K^{n-1}P(S_m \leq t) \\
& \leq M (b_{n+1}+b_{n})K^n P(S_{k_n} \leq t) =M \, b_{n+2} \, K^n P(S_{k_{n+1}} \leq t)
\end{align*}
This finishes the proof of (\ref{induc}).

From (\ref{induc}) we infer that \eqref{fn-bounded-an} holds with $a_0=a_1=1$ and  $$a_n=b_{n+1}K^{n-1}P(S_{k_n} \leq T) \quad \mbox{for} \quad n \geq 2.$$
The fact that $\sum_{n \geq 0}a_n^{1/p}<\infty$ for any $p>0$ follows immediately, since $$b_n=\frac1{\sqrt{5}}\left[\left(\frac{1+\sqrt{5}}2\right)^{n}-\left(\frac{1-\sqrt{5}}2\right)^{n}\right]$$
and by Lemma 17 of \cite{dalang99}, we know that for any
$a>1$ and for any $p>0$,
$$\sum_{n=1}^{\infty}a^{n/p} P(S_n\le T)^{1/p}<\infty.$$
\qed

\appendix

\section{Peszat's condition}
\label{appendix-peszat-cond}

In this section, we show that condition \eqref{eq:51} fails when $\mu$ is given by \eqref{def-mu}.

\begin{lemma}
For any $H<1/2$,
$$\sup_{\eta \in \bR}\int_{\bR}\frac{1}{1+|\xi-\eta|^2}|\xi|^{1-2H}\,d\xi=\infty.$$
\end{lemma}

\noindent {\bf Proof:} We have
\begin{align*}
& \sup_{\eta \in \bR}\int_{\bR}\frac{1}{1+|\xi-\eta|^2}\,|\xi|^{1-2H}d\xi=
\sup_{\eta \in \bR}\int_{\bR}\frac{1}{1+|\xi|^2}\,|\xi+\eta|^{1-2H}d\xi \\
& \quad \geq \sup_{\eta \geq 0}\int_{0}^{\infty}\frac{1}{1+|\xi|^2}\,(\xi+\eta)^{1-2H}d\xi=\lim_{\eta \to \infty}\int_{0}^{\infty}\frac{1}{1+|\xi|^2}\,(\xi+\eta)^{1-2H}d\xi,
\end{align*}
since the last integral is an increasing function of $\eta$. The conclusion follows by the monotone convergence theorem, since $\lim_{\eta \to \infty}(\xi+\eta)^{1-2H}=\infty$. $\Box$

\section{Martingale Measure}
\label{appendix-mart}

The following definition was introduced in \cite{walsh86} for real-valued processes. Here we extend it to complex-valued processes.

\begin{definition}
{\rm A complex-valued process $\{M_t(A);t \geq 0,A \in \cB_b(\bR)\}$ is a {\em martingale measure} with respect to a filtration $\{\cF_t\}_{t \geq 0}$ if\\
(a) for any $A \in \cB_b(\bR)$, $\{M_t(A)\}_{t \geq 0}$ is a square-integrable complex-valued martingale with respect to $\{\cF_t\}_{t \geq 0}$, with $M_0(A)=0$ a.s.; \\
(b) for any $t>0$, $\{M_t(A);A \in \cB_{b}(\bR)\}$ is a $\sigma$-finite $L^2(\Omega)$-valued signed measure, in the sense that it satisfies the following two properties:\\
(b1) $M_t(A \cup B)=M_t(A)+M_t(B)$ a.s. for any disjoint sets $A,B \in \cB_b(\bR)$;\\
(b2) there exists a sequence $(E_k)_{k} \subset \cB_b(\bR)$ with $E_k \subset E_{k+1}$ for all $k$ and $\cup_{k}E_k=\bR$, such that for any $k$,
$\sup_{A \in \cB_k}E|M_t(A)|^2<\infty$ and $E|M_t(A_n)|^2 \to 0$ for any sequence $(A_n)_n \subset \cB_k$ with $A_n \downarrow \emptyset$, where $\cB_k=\{B \in \cB(\bR);B \subset E_k\}$.
}
\end{definition}

\begin{lemma}
The process $M=\{M_t(A);t \geq 0,A \in \cB_b(\bR)\}$ defined by \eqref{eq:333} is a complex-valued martingale measure with respect to $\{\cF_t\}_{t \geq 0}$ with covariation
$$\langle M(A),M(B)\rangle_{t}=t\mu(A \cap B).$$
\end{lemma}

\noindent {\bf Proof:} By approximation with simple functions, we can define the integral $\cM(\varphi):=\int_{\bR^2}\varphi(\tau,\xi)\cM(d\tau,d\xi)$ as an element in $L_{\bC}^2(\Omega)$, for any $\varphi \in L_{\bC}^2(\bR^2,\Pi)$. This integral has the property that for any $\varphi,\psi \in L_{\bC}^2(\bR^2,\Pi)$,
$$E[\cM(\varphi) \overline{\cM(\psi)}]=\frac{1}{2\pi}\int_{\bR^2} \varphi(\tau,\xi) \overline{\psi(\tau,\xi)}d\tau \mu(d\xi).$$

Therefore, for any $s<t$ and $A \in \cB_b(\bR)$, $M_t(A)-M_s(A)=\cM(\cF 1_{(s,t]} \cdot 1_{A})$ is orthogonal
to $X_u(\varphi)=\cM(\cF 1_{[0,u]} \cdot \cF \varphi)$, for any $u<s$ and $\varphi \in \cD(\bR)$. Since $\{\cM(\varphi); \varphi \in L_{\bC}^2(\bR^2,\Pi)\}$ is a Gaussian process, it follows that $M_t(A)-M_s(A)$ is independent of $\cF_s$. In addition, $E[M_t(A)]=0$ for any $t>0$. Hence,
$\{M_t(A)\}_{t \geq 0}$ is a martingale with respect to $\{\cF_t\}_{t \geq 0}$. This proves that $M$ satisfies property (a). Property (b) follows immediately from
\eqref{cov-M}. The statement about the covariation follows from Proposition 4.1 of \cite{BGP12}. \qed

\section{Some auxiliary results}
\label{appendix-elem}

\begin{lemma}
\label{elem-lem1}
The integral $\int_{0}^{\infty}(1-\cos x)x^{-\alpha-1}dx=:c_{\alpha}$ converges if and only if $\alpha \in (0,2)$. In this case,
$$c_{\alpha}=
\left\{
\begin{array}{ll}
\alpha^{-1}\Gamma(1-\alpha) \cos(\pi \alpha/2) & \mbox{if $\alpha \in (0,1)$} \\
-\alpha^{-1}(\alpha-1)^{-1}\Gamma(2-\alpha)\cos(\pi \alpha/2) & \mbox{if $\alpha \in (1,2)$}\\
\pi/2 & \mbox{if $\alpha=1$}.
\end{array} \right.
$$
For any $\alpha \in (0,2)$ and for any $\xi>0$, we have
\begin{equation}
\label{cos-int}
\int_{0}^{\infty}\frac{1-\cos (\xi x)}{x^{\alpha+1}}dx=c_{\alpha}\xi^{\alpha}.
\end{equation}
\end{lemma}

\noindent{\bf Proof:} To see that the integral converges if $\alpha \in (0,2)$, we use the fact that
$1-\cos x \leq x^2/2$ when $|x| \leq 1$, and $1-\cos x \leq 2$ when $|x|>1$.

For the other implication, note that $\sin^2 r \geq r^2 \sin^2 1$ for any $r \in [0,1]$. Hence, for any $x \in [0,2]$, $1-\cos x =2 \sin^2(x/2) \geq  2^{-1}(\sin^2 1)x^2$ and
$$\infty>\int_{0}^{\infty}\frac{1-\cos x}{x^{\alpha+1}}dx \geq 2^{-1}(\sin^2 1) \int_0^2 x^{2-\alpha-1}dx \quad
\mbox{implies that} \ \alpha<2.$$ It remains to show that the
integral diverges if $\alpha \leq 0$. In fact,
applying a change of variables, we have
\[
\int_{0}^{\infty}\frac{1-\cos x}{x^{\alpha+1}}dx =
\frac{2}{\pi^{\alpha}}
\int_{0}^{\infty}\frac{\sin^2(x\pi/2)}{x^{\alpha+1}}dx.
\]
Now, observe that for some small $\varepsilon>0$ we will have that
\begin{align*}
\int_{0}^{\infty}\frac{\sin^2(x\pi/2)}{x^{\alpha+1}} \, dx & \geq
\sum_{n \geq 1\, n \, {\rm odd}} \int_{n-\varepsilon}^{n+\varepsilon}
\frac{\sin^2(x\pi/2)}{x^{\alpha+1}}\, dx \\
& \geq 2 \varepsilon \sin^2(\pi(\varepsilon+1)/2) \sum_{n\geq 1\, n\, {\rm odd}}
\frac{1}{(n+\varepsilon)^{\alpha+1}},
\end{align*}
and the latter series diverges whenever $\alpha\leq 0$. Here we used the fact that
if $x \in (n-\varepsilon,n+\varepsilon)$ and $n$ is an odd integer, then $\sin^2 (x\pi/2) \geq \sin^2 (\pi(n+\varepsilon)/2)=\sin^2(\pi(1+\varepsilon)/2)$.

In order to deduce the explicit formulas for $c_\alpha$, we first
integrate by parts, yielding
$$\int_{0}^{\infty}\frac{1-\cos y}{y^{\alpha+1}}dy=\frac{1}{\alpha}\int_{0}^{\infty}\frac{\sin y}{y^{\alpha}}dy.$$

By formula 3.761-4 of \cite{gradshteyn-ryzhik07}, we know that
\begin{equation}
\label{sin-integral}
\int_0^\infty x^{\mu-1}\,{\sin x}\,
dx=\frac{\pi}{2\,\Gamma(1-\mu)\,\cos(\mu\pi/2)}, \quad \mbox{for} \ \mu\in
(0,1).
\end{equation}

If $\alpha \in (0,1)$, using \eqref{sin-integral} with $\mu=1-\alpha$, we have
$$\int_{0}^{\infty}\frac{\sin y}{y^{\alpha}}dy =
\frac{\pi}{2\Gamma(\alpha)\sin(\pi \alpha/2)}=\Gamma(1-\alpha)\cos\left(\frac{\pi \alpha}{2}\right).$$
where for the second equality we used the identity:
$$\Gamma(\alpha)\,\Gamma(1-\alpha)=\frac{\pi}{\sin(\pi \alpha)} \quad \mbox{for any} \ \alpha \in (0,1)$$
(see formula (6.1.17) of \cite{avramowitz-stegun64}).

The calculation of $c_{\alpha}$ in the case $\alpha \in [1,2)$ is
given on pages 568-569 of \cite{feller}.

Eventually, relation \eqref{cos-int} follows using the change of
variables $y=\xi x$. \qed

\begin{lemma}
\label{elem-lem2}
For any $0<H<1/2$ and $\xi \in \bR$, we have:
$$\int_{\bR}\frac{|1-e^{-i\xi x}|^2}{|x|^{2-2H}}dx=|\xi|^{1-2H}\frac{2\,\Gamma(2H+1)\sin(\pi H)}{H(1-2H)}.$$
\end{lemma}

\noindent {\bf Proof:} Assume that $\xi>0$. (The case $\xi<0$ is similar.)
$$\int_{\bR}\dfrac{|1-e^{-i\xi x}\,|^2}{|x|^{2-2H}}dx=4\int_0^\infty
\frac{1-\cos (\xi
x)}{x^{2-2H}}dx=\xi^{1-2H}\frac{4\Gamma(2H)}{1-2H}\cos(\pi
(1-2H)/2),$$ where for the second equality we used Lemma
\ref{elem-lem1} with $\alpha=1-2H$. The conclusion follows since
$\Gamma(2H+1)=2H\Gamma(2H)$. \qed

\end{document}